\documentclass[letterpaper,10pt]{IEEEtran}

\usepackage{graphics}
\usepackage{epsfig}
\usepackage{amsmath}
\usepackage{amssymb}
\usepackage{cite}
\usepackage{color}
\usepackage{enumerate}
\usepackage{graphicx}
\usepackage{mathtools}
\usepackage{amsthm}
\usepackage{cancel}
\usepackage{url}
\usepackage{derivative}

\newcommand{\setmap}[3]{#1:#2 \mathrel{\vcenter{\mathsurround0pt
\ialign{##\crcr
		\noalign{\nointerlineskip}$\rightarrow$\crcr
		\noalign{\nointerlineskip}$\rightarrow$\crcr
		}}}%
		#3}

\newcommand{\naturals}{\mathbb{N}}
\newcommand{\real}{\mathbb{R}}
\newcommand{\R}{\mathbb{R}}
\newcommand{\realnonneg}{\mathbb{R}_{\ge 0}}
\newcommand{\realpos}{\mathbb{R}_{> 0}}
\newcommand{\until}[1]{[#1]}
\newcommand{\map}[3]{#1:#2 \rightarrow #3}
\newcommand{\qedA}{~\hfill \ensuremath{\square}}

\newcommand{\interior}{\operatorname{int}}

\newcommand{\longthmtitle}[1]{\mbox{}{\textit{(#1):}}}

\newcommand{\setdefb}[2]{\big\{#1 \; | \; #2\big\}}
\newcommand{\setdefB}[2]{\Big\{#1 \; | \; #2\Big\}}
\newcommand{\Ball}[2]{\mathfrak B_{#1}(#2)}
\newcommand*{\SetSuchThat}[1][]{}
\newcommand*{\MvertSets}{%
    \renewcommand*\SetSuchThat[1][]{%
        \mathclose{}%
        \nonscript\;##1\vert\penalty\relpenalty\nonscript\;%
        \mathopen{}%
    }%
}
\MvertSets
\DeclarePairedDelimiterX \Set [2] {\lbrace}{\rbrace}
    {\,#1\SetSuchThat[\delimsize]#2\,}

\newcommand{\Ac}{\mathcal{A}}
\newcommand{\Cc}{\mathcal{C}}
\newcommand{\Dc}{\mathcal{D}}

\newcommand{\Kc}{\mathcal{K}}

\newcommand{\Pc}{\mathcal{P}}
\newcommand{\Uc}{\mathcal{U}}
\newcommand{\Sc}{\mathcal{S}}

\newcommand{\Wc}{\mathcal{W}}

\newcommand{\Ic}{\mathcal{I}}

\newtheorem{theorem}{Theorem}[section]
\newtheorem{lemma}[theorem]{Lemma}
\newtheorem{corollary}[theorem]{Corollary}
\newtheorem{proposition}[theorem]{Proposition}

\theoremstyle{definition}
\newtheorem{definition}{Definition}

\newtheorem{remark}{Remark}
\newtheorem{example}{Example}

\newtheorem{problem}{Problem}

\newcommand{\ba}{\mathbf{a}}
\newcommand{\bb}{\mathbf{b}}

\newcommand{\be}{\mathbf{e}}
\renewcommand{\bf}{\mathbf{f}}
\newcommand{\bg}{\mathbf{g}}

\newcommand{\bk}{\mathbf{k}}

\newcommand{\bu}{\mathbf{u}}

\newcommand{\bw}{\mathbf{w}}
\newcommand{\bx}{\mathbf{x}}

\newcommand{\bB}{\mathbf{B}}

\newcommand{\bG}{\mathbf{G}}

\newcommand{\bI}{\mathbf{I}}

\newcommand{\bT}{\mathbf{T}}

\newcommand{\bmu}{\boldsymbol{\mu}}

\newcommand{\bzero}{\mathbf{0}}

\newcommand{\des}{{\operatorname{d}}}

\newcommand{\T}{^\top}

\newcommand{\son}{{\operatorname{stg}}}
\newcommand{\tri}{{\operatorname{tri}}}
\newcommand{\mir}{{\operatorname{m}}}

\newcommand{\clf}{{\operatorname{CLF}}}
\newcommand{\cbf}{{\operatorname{CBF}}}

\newcommand{\gau}{{\operatorname{gau}}}
\newcommand{\Lie}{\mathcal{L}}

\newcommand{\defeq}{\triangleq}

\allowdisplaybreaks

\title{\LARGE  Construction of Universal Safe Stabilizing Formulas \\ with Weighted Centroid Methodology
  \thanks{This work was supported by AFOSR Award FA9550-23-1-0740.}  }

\author{Pio Ong \quad Max H. Cohen \quad Aaron D. Ames \quad Jorge Cort\'es%
\thanks{Pio Ong is with the Amazon SAF Lab under Industrial Robotics Group, Amazon.com, Inc., Pasadena, CA 91125, USA; \texttt{pioong@amazon.com}. This work is unrelated to his position at Amazon. The work was initiated while he was with the University of California, San Diego, and it was developed primarily while he was with the California Institute of Technology.}
\thanks{Max H. Cohen is with the Department of Electrical and Computer Engineering, North Carolina State University, Raleigh, NC 27695, USA; \texttt{mhcohen2@ncsu.edu}}%
\thanks{Aaron D. Ames is with the Department of Mechanical and Civil Engineering, California Institute of Technology, Pasadena, CA 91125, USA; \texttt{ames@caltech.edu}}%
  \thanks{Jorge Cort\'es is with Department of Mechanical
    and Aerospace Engineering, University of California, San Diego, CA
    92093, USA; \texttt{cortes@ucsd.edu}}%
}

\begin{document}
\maketitle

\begin{abstract}
    Control Lyapunov and control barrier functions  ensure stability and safety, respectively, by imposing state-dependent affine constraints on the control input. Controller design can therefore be viewed as the problem of selecting an admissible input from the set defined by these constraints. This paper adopts this perspective and develops a weighted-centroid methodology for constructing explicit feedback formulas that are provably smooth. The selection method relies on the fact that weighted centroids of a convex set lie within the set. We establish conditions on the weight functions that guarantee smoothness and continuity for several classes of constraints, including single constraints, mirrored constraints, linearly independent multiple constraints, and two-constraint Gaussian constructions. The framework recovers Sontag's universal formula and generates other new explicit formulas. When specialized to control-theoretic settings, the resulting formulas yield smooth safety filters, stabilizing controllers, and safe-stabilizing controllers. Simulation and hardware experiments on a small quadrotor illustrate the practical benefits of the proposed construction.
\end{abstract}

\section{Introduction}
Certificate-based methods provide powerful tools for enforcing and certifying desirable closed-loop dynamical properties, such as stability and safety. Control Lyapunov functions (CLFs)~\cite{EDS:83,RAF-PVK:96} certify the existence of stabilizing inputs, while control barrier functions (CBFs)~\cite{ADA-XX-JWG-PT:17,ADA-SC-ME-GN-KS-PT:19} certify the existence of inputs that render a prescribed safe set forward invariant. In both cases, these certificates translate the desired closed-loop properties into state-dependent constraints on the control inputs. As a result, controller design can be viewed as a \textit{selection} problem consisting of selecting at each state an admissible control input from the set defined by the constraints.

These selection problems can be approached using optimization, formulating the control design as a state-wise optimization subject to the control-theoretic constraints. Such ideas have long been employed in nonlinear control~\cite{RAF-PVK:96}, but have become especially prominent over the past decade in the safety-critical control literature via CLF-CBF-based quadratic programs~\cite{ADA-XX-JWG-PT:17}. This optimization framework is appealing due to its modularity, computational tractability, and its ability to accommodate multiple constraints, making it well suited for real-time safety-critical applications. Despite its advantages, optimization-based controllers often fail to be differentiable~\cite{PM-AA-JC:25,BJM-MJP-ADA:15,DRA-HL-DP:25}.  This nonsmoothness is already evident from the closed-form expressions in the single-constraint case~\cite{MJ:18}, and persists in more general explicit CBF-based safety filters~\cite{PM-SSM-PO-LY-ED-JWB-ADA:26-ifacwc}.

In many settings, it is desirable for the feedback controller to be smooth. In particular, smoothness enables recursive and compositional control constructions in which the controller itself is an intermediate design step within a broader control design framework or hierarchy. Examples include backstepping \cite{MK-IK-PK:95}, which has recently been extended to safety-critical control \cite{AJT-PO-TGM-ADA:22-cdc,MHC-PO-ADA:26-acc}. Other examples include the backup CBF framework \cite{TG-MM-AS-PN-EF-ADA:20} --- in which a safe backup policy is used to dynamically generate an implicit control invariant set --- where the backup controller must be continuously differentiable to reliably compute gradients of the resulting implicit CBF.

Smooth feedback design has long appeared in the literature of nonlinear stabilization. In particular, Artstein's theorem~\cite{ZA:83} establishes the existence of a continuous stabilizing controller from a CLF. Although not stated explicitly, the same selection method can be extended more generally to obtain smooth feedback~\cite{PO-BC-LS-JC:23-auto}.
Unfortunately, these existence results are nonconstructive. Sontag's universal formula~\cite{EDS:89a} addresses this limitation by providing an explicit feedback construction from a CLF. The formula, however, is tailored to a single CLF constraint. Subsequent works \cite{YL-EDS:91} and \cite{YL-EDS:95} extend the universal formula to handle additional constraints such as $u\in[0,\infty)$ and $u\in[0,1]$, respectively. However, these methods do not handle additional constraints that are state-dependent, such as those induced by CBFs.

Motivated in part by the introduction of CBFs, more recent works have revisited universal-formula constructions. Earlier work applies the Sontag's formula in settings where both stability and safety objectives can be encoded together with control Lyapunov barrier functions~\cite{MZR-BJ:16,PB-CMK:20}.
The implicit-function-theorem viewpoint underlying Sontag's formula has been used to characterize and construct smooth safety filters for single affine constraints~\cite{MHC-PO-GB-ADA:23-csl}. Other works~\cite{ML-ZS:23,ML-ZS-PJWK-SW:24} study how the Sontag's formula can be tuned toward minimum-norm controllers while retaining smoothness or inverse optimality~\cite{MK:24}. Extensions to safe stabilization with both CLF and CBF constraints have also been developed. For instance, \cite{BW-MK:26} proposes a Sontag-type formula that handles both constraints, albeit for single-input systems. Alternatively, multiple constraints with multiple inputs can be handled by reformulating the selection problem as an unconstrained optimization problem whose optimizer varies smoothly with the state~\cite{PM-JC-EDS:26}. While this provides a smooth selection mechanism, it does not yield an explicit closed-form formula.

\subsubsection*{Statement of Contributions}
The present paper develops a general weighted-centroid methodology for constructing explicit smooth feedback formulas that select admissible control inputs from state-dependent affine constraints. The work builds on our earlier conference version of the paper~\cite{PO-JC:19-cdc}, which introduced a Gaussian-weighted centroid formula for smooth safe stabilization with one CLF and one CBF. Here, we generalize the core idea into a broader framework in which the choice of weight function becomes a mechanism for generating explicit feedback formulas.

The first contribution is the introduction of the weighted-centroid perspective as a constructive methodology for control design.
In contrast to~\cite{PO-JC:19-cdc}, where smoothness followed from known closed-form expressions for truncated Gaussian means, the present smoothness results apply to general classes of weight functions.
In particular, we identify sufficient conditions for smoothness of the admissible control maps and weight functions, covering the following cases: single affine constraints, mirrored constraints (two affine constraints with anti-parallel normals), multiple constraints with linearly independent constraint normals, structured linearly dependent constraints arising from mirrored pairs, and the general two affine constraints using Gaussian weights. For linearly independent cases, we provide continuity guarantees at states where constraint normals vanish, provided the remaining nonzero constraint normals remain linearly independent.

The second contribution is the use of the proposed methodology to derive explicit feedback formulas, beyond Gaussian-weighted constructions in~\cite{PO-JC:19-cdc}. Notably, we show that Sontag's universal formula can be recovered as a weighted centroid of the admissible control set induced by a single affine constraint. This interpretation allows us to extend Sontag's formula to handle mirrored constraints while retaining an explicit closed-form expression, rather than an optimization-based formulation. Beyond Sontag-type constructions, we derive additional explicit formulas using Gaussian weights and various polynomial weights. We also develop a compactly supported formula that respects state-independent input bounds.

The third contribution is the specialization of the framework to control-theoretic settings with CLFs and CBFs. The resulting constructions yield explicit smooth formulas for safety filtering, stabilization, and safe stabilization by satisfying the corresponding CLF and/or CBF constraints under suitable compatibility assumptions. For the CLF stabilization setting,
we characterize classes of weight functions for which the resulting formula is continuous at the origin under the small control property, and propose a state-dependent scaling method for modifying existing weight functions to fit this class.

Finally, we demonstrate the practical implications of the proposed methodology through both simulation and hardware experiments on a quadrotor. The simulation study uses higher-order tracking dynamics to illustrate how nonsmooth safety-filter commands can repeatedly excite tracking transients, whereas the proposed smooth formula avoids this effect. Complementary hardware experiments on a small quadrotor show that our approach provides an inherent robustness margin against residual tracking and implementation errors.

\section{Admissible Control Map Selection}
\label{sec:problem}
This section introduces the problem of designing a controller to satisfy state-dependent control constraints. Consider a nonlinear control-affine system\footnote{We use $\naturals$, $\real$, and $\realnonneg$, to denote the sets of natural, real, and nonnegative real numbers, respectively. The space $\real^n$ is the $n$-dimensional Euclidean space, and $\|\cdot\|$ denotes the Euclidean norm. For $p\in\naturals$, let $\until p$ denote the set of consecutive numbers $\{1,\dots, p\}$. Given $\bx\in\real^n$ and $\delta>0$, $\Ball{\delta}{\bx}$ is the open ball of radius $\delta$ centered at $\bx$. The zero vector is denoted by $\bzero$, with dimension inferred from context.
For a set $\Sc$, we denote its interior, boundary and complement by $\interior(\Sc)$, $\partial\Sc$, and $\Sc^c$, respectively. A set-valued map $\setmap{\Uc}{\real^n}{\real^m}$ is a mapping of $\bx\in\real^n$ to a set $\Uc(\bx)\subseteq \real^m$.
For a continuously differentiable function $\map{h}{\real^n}{\real}$, we denote its Lie derivative along a collection of vector fields $\map{\bg}{\real^n}{\real^{n\times m}}$ by $\Lie_\bg h(\bx) \defeq \frac{\partial h}{\partial\bx}(\bx)\bg(\bx)$.
In this paper, all integrals are understood in the Lebesgue sense, and a function is called integrable if it is absolutely integrable.
A function $\map{\alpha}{\realnonneg}{\realnonneg}$ is of class-$\Kc$ if $\alpha(0)=0$ and $\alpha$ is strictly increasing. Furthermore, $\alpha$ is a class-$\Kc_\infty$ function if it is of class-$\Kc$ and ${\lim_{r\to\infty}\alpha(r) = \infty}$. An extended class-$\Kc$, denoted by $\Kc^e$, function $\alpha$ has the property of a class-$\Kc$ function while has its domain and range extended to include $0$ in their interior.}:
\begin{equation}
\label{sys:ctrl-affine}
\dot \bx = \bf(\bx) +\bg(\bx)\bu,
\end{equation}
where $\bx\in\real^n$ is the state and $\bu\in\real^m$ is the control input. The drift vector field $\map{\bf}{\real^n}{\real^n}$ and the control matrix-valued function $\map{\bg}{\real^n}{\real^{n\times m}}$ are assumed to be smooth. Given a locally Lipschitz feedback controller $\map{\bk}{\real^n}{\real^m}$, the closed-loop system is:
$\dot{\bx} = \bf(\bx) + \bg(\bx)\bk(\bx)$.

In many situations, the choice of $\bu$ is limited to a specific set ${\Uc \subset \real^m}$. In the simplest case, the set $\Uc$ is fixed, e.g., when the input has hard actuation limits. In general, however, the constraint set is described by a set-valued map $\setmap{\Uc}{\real^n}{\real^m}$ that captures dependency on the state. This naturally arises in control design problems that enforce performance or safety specifications. For example, control Lyapunov functions (CLFs) \cite{EDS:83} and control barrier functions (CBFs) \cite{ADA-XX-JWG-PT:17} impose inequality constraints on $\bu$ that depend explicitly on the state.
Although the results in this paper apply to general state-dependent affine constraints, we focus on those induced by CLFs and CBFs, given their relevance in control theory.

\subsection{Control Lyapunov and Control Barrier Functions}
The concept of a CLF builds on Lyapunov theory by incorporating the control input~$\bu$ to certify the existence of stabilizing controllers. Throughout this section, we assume that the system \eqref{sys:ctrl-affine} admits an equilibrium under zero input. Without loss of generality, we shift coordinates so that the equilibrium of interest is located at the origin $\bx=0$, i.e., $\bf(0)=0$.

\begin{definition}\longthmtitle{Control Lyapunov Function \cite[Section 3.3.1]{RAF-PVK:96}}
\label{def:clf}
A continuously differentiable function  ${\map{V}{\real^n}{\real}}$ is a \emph{control Lyapunov function} for the control-affine system \eqref{sys:ctrl-affine} if there exists class-$\Kc_{\infty}$ functions $\underline{\alpha}$ and $\overline{\alpha}$ such that for all $\bx\in\real^n$:
$$
\underline{\alpha}(\|\bx\|)\leq V(\bx) \leq \overline{\alpha}(\|\bx\|),
$$
and, there exists a class-$\Kc$ function $\alpha$ such that, for each $\bx\in\real^n\setminus \{0\}$, there exists $\bu\in\real^m$ satisfying \emph{CLF constraint}:
\begin{equation}\label{eq:clf_deriv}
\Lie_\bf V(\bx)+\Lie_\bg V(\bx)\bu< -\alpha(V(\bx))
\end{equation}
\end{definition}

At each $\bx\in\real^n$, the CLF constraint~\eqref{eq:clf_deriv} ensures the existence of a control $\bu$ that makes the time derivative of the function $V$ negative, as required to make $V$ a Lyapunov function.
 Importantly, CLFs facilitate the selection of control inputs to enforce stability of the closed-loop system.

\begin{lemma}\longthmtitle{Stabilizing Controller}\label{lem:stab}
Consider the control-affine system~\eqref{sys:ctrl-affine} and a CLF~$V:\real^n \rightarrow \real$. Define a state-dependent admissible control map $\setmap{\Omega_\text{CLF}}{\real^n}{\real^m}$ such that for all $\bx\neq 0$:
$$
\Omega_\clf(\bx) = \setdefb{\bu\in\real^m}{\Lie_\bf V(\bx)+\Lie_\bg V(\bx)\bu\leq -\alpha(V(\bx))}.
$$
Given a function $\map{\bk}{\real^n}{\real^m}$ such that $\bk(\bzero)=\bzero$, $\bk$ is locally Lipschitz on $\real^n\setminus\{\bzero\}$, and $\bk(\bx)\in\Omega_\clf(\bx)$ for each $\bx\in \real^n \setminus \{\bzero\}$, the feedback law $\bu=\bk(\bx)$ renders the origin asymptotically stable.~\hfill$\blacksquare$
\end{lemma}

The set-valued map $\Omega_\clf$ specifies, at each $\bx$, the set of all admissible inputs $\bu$ that lead to the decrease of~$V$. Note that although the constraint~\eqref{eq:clf_deriv} is stated with a strict inequality, the admissible control set $\Omega_{\clf}$ is defined using
the non-strict condition ($\leq$). The strict inequality in the definition not only guarantees that the
set $\Omega_{\clf}(\bx)$ is nonempty but also ensures the existence of a smooth feedback controller satisfying the constraint, see also Lemma~\ref{lem:artstein} discussed shortly.

We next review the concepts of safety and CBFs. For \eqref{sys:ctrl-affine}, suppose $\Sc\subset \real^n$ is a safe set and the control objective is to design a feedback controller such that the safe set $\Sc$ is forward invariant for the closed-loop system. In other words, the set $\Sc$ is \emph{safe} if the solution to the closed-loop system is such that $\bx(t)\in\Sc$ at all time in its interval of existence for any given initial condition $\bx_0\in\Sc$. The concept of a CBF facilitates the construction of controllers enforcing safety.

\begin{definition}\longthmtitle{Control Barrier Function \cite{ADA-XX-JWG-PT:17}}
\label{def:cbf}
{\rm
A continuously differentiable function $\map{h}{\real^n}{\real}$ is a \emph{control barrier function} for the set $\Sc$ under the control-affine system~\eqref{sys:ctrl-affine} if its 0-superlevel set defines $\Sc$:
\begin{equation}
\label{eq:safeset}
\begin{aligned}
\Sc = & \setdefb{\bx\in\real^n}{h(\bx)\geq 0};
\end{aligned}
\end{equation}
and,  there exists a class-$\Kc^e$ function $\alpha_{h}$ such that, for each ${\bx\in\Sc}$, there exists ${\bu\in\real^m}$ satisfying the \emph{CBF constraint}:
\begin{equation}\label{eq:cbf_deriv}
\Lie_\bf h(\bx)+\Lie_\bg h(\bx)\bu > -\alpha_{h}(h(\bx)),
\end{equation}
}
\end{definition}

The concept of a CBF builds on that of a barrier certificate \cite{SP-AJ:04}, which has traditionally been used for safety verification of autonomous systems, by extending such a certificate function to incorporate control inputs.
The CBF constraint~\eqref{eq:cbf_deriv} guarantees that there exists a control $\bu$ that does just that. Furthermore, the function $\alpha_{h}$ on the right-hand side of \eqref{eq:cbf_deriv} allows the value of $h$ to decrease at points where $h(\bx)>0$, i.e., in the interior of~$\Sc$. Consequently, there are more possible choices of control input $\bu$ that help render the set $\Sc$ safe.
We formally state the safety result using CBFs below.

\begin{lemma}\longthmtitle{Safeguarding Controller}\label{lem:safe}
Consider the control-affine system~\eqref{sys:ctrl-affine} and a CBF $h:\real^n \rightarrow \real$. Define a state-dependent admissible control map $\setmap{\Omega_\cbf}{\Sc}{\real^m}$ as:
$$
\Omega_\cbf(\bx) = \setdefb{\bu\in\real^m}{\Lie_\bf h(\bx)+\Lie_\bg h(\bx)\bu \geq -\alpha_{h}(h(\bx))}.
$$
Given a function  $\map{\bk}{\real^n}{\real^m}$ such that $\bk$ is locally Lipschitz on $\Sc$ and $\bk(\bx)\in\Omega_\cbf(\bx)$ for each $\bx\in \Sc$, the feedback $\bu=\bk(\bx)$ renders the set $\Sc$ safe.~\hfill$\blacksquare$
\end{lemma}
Lemma~\ref{lem:safe} uses a CBF to guarantee safety by limiting the choice of inputs at each state. Similar to CLFs, CBFs induce a set-valued map $\Omega_\cbf$ that specifies admissible control inputs. The use of a strict inequality in the CBF constraint~\eqref{eq:cbf_deriv}, while the safety result only requires the non-strict condition defining $\Omega_\cbf$, parallels the discussion following Lemma~\ref{lem:stab}.

\subsection{Problem Statement}
Motivated by the CLF and CBF constructions above, we consider a general class of affine inequality constraints. Abstracting this structure, we define  \textit{admissible control maps}
$\setmap{\Omega_i}{\Sc_i}{\real^m}$ on domains
$\Sc_i\subseteq\real^n$ as:
\begin{equation}\label{eq:ctrl-affine-feas}
\Omega_i(\bx) = \setdefb{\bu \in \real^m}{a_i(\bx)+\bb_i(\bx)\bu\leq 0~(\text{or}<0)}
\end{equation}
with state-dependent functions $\map{a_i}{\real^n}{\real}$ (\textit{constraint bias}) and  $\map{\bb_i}{\real^n}{\real^{1\times m}}$ (\textit{constraint normal}).
Given a collection of $q$ admissible control maps~$\{\Omega_i\}_{i=1}^q$, we address the \emph{selection} problem: finding a function $\map{\bk}{\real^n}{\real^m}$ such that
$\bk(\bx)\in\Omega_i(\bx)$, $\forall i\in\until{q}$. To this end, we introduce \textit{compatibility}.
\begin{definition}\longthmtitle{Compatibility}\label{def:compatibility}
{\rm
The admissible control maps $\{\Omega_i\}_{i=1}^q$ are \emph{(strictly) compatible} on $\Sc$ if, for each $\bx\in\Sc$, there exists a control input $\bu\in\real^m$ such that $a_i(\bx)+\bb_i(\bx)\bu \leq 0~(\text{or} <0)$ for all $i\in\until{q}$.~\hfill$\diamond$
}
\end{definition}
Indeed, compatibility among the admissible control maps is required for the existence of a selection function $\bk$. However, we will make a slightly stronger assumption of strict compatibility among the admissible control maps, which is closely related to the possibility of finding a selection function $\bk$ that is also smooth. The existence of a smooth selection is guaranteed by the following result.

\begin{lemma}\longthmtitle{Generalization of Artstein's Theorem~\cite{ZA:83}}\label{lem:artstein}
Let $\{\setmap{\Uc_i}{\real^n}{\real^m}\}_{i=1}^q$ be a set of $q$ admissible control maps:
\begin{equation}\label{eq:feas}
\Uc_i(\bx) = \setdefb{\bu \in \real^m}{\gamma_i(\bx,\bu)\leq 0~(\textup{or}<0)}
\end{equation}
with constraint functions $\map{\gamma_i}{\real^n \times \real^m}{\real}$. Consider the intersection $\Uc = \cap_{i=1}^q \Uc_i $ and assume the following hold:
\begin{enumerate}
\item $\gamma_i$ are continuous for all $i\in\until q$;
\item $\Uc$ is convex-valued;
\item the maps $\{\Uc_i\}_{i=1}^q$ are strictly compatible on $\Sc = \cap_{i=1}^q \Sc_i$
\end{enumerate}
Then there exists a $\mathcal C^\infty$ selection function $\map{\bk}{\Sc}{\real^m}$ such that $\bk(\bx)\in\Uc(\bx)$ for all $\bx\in\Sc$.~\hfill $\qedA$
\end{lemma}

In the context of control, where the selection function $\bk$ is the controller, it is important for $\bk$ to be smooth or, at the very least, locally Lipschitz continuous, to guarantee existence and uniqueness of solutions. Similarly, many results in control theory often require Lipschitzness, e.g., Lemmas~\ref{lem:stab} and~\ref{lem:safe}, and others may even require smoothness, e.g., backstepping~\cite{MK-IK-PK:95}. We summarize the design problem below.

\begin{problem}\longthmtitle{Smooth Selection Problem}
\label{prob:selection}
Consider the control-affine system~\eqref{sys:ctrl-affine} with admissible control maps $\{\Omega_i\}_{i=1}^q$ defined in~\eqref{eq:ctrl-affine-feas}. Assume that the admissible control maps are strictly compatible on $\Sc = \cap_{i=1}^q \Sc_i$. Design a smooth feedback controller $\map{\bk}{\Sc}{\real^m}$ such that $\bk(\bx)\in \Uc(\bx) =\cap_{i=1}^q \Omega_i(\bx)$ for all $\bx\in \Sc$.~\hfill$\bullet$
\end{problem}

The remainder of the paper is organized as follows. Section~\ref{sec:method} introduces our proposed selection methodology and develops the construction in the case of a single constraint or single-input systems. Section~\ref{sec:general-case} extends the approach to the general case with multiple constraints and inputs. Section~\ref{sec:control} specializes the results to control applications, providing corollaries for the simultaneous enforcement of stabilization and safety objectives. Finally, Section~\ref{sec:experiments} demonstrates the proposed approach through a higher-order quadrotor simulation and hardware experiments on a small quadrotor platform.

\section{Weighted Centroid Selection Methodology}
\label{sec:method}
Our approach to the selection problem hinges on the convexity (in~$\bu$) of the admissible control maps. In particular, we rely on the fact that a weighted centroid of a convex set is always inside the set \cite[p. 25]{SB-LV:09}. This property suggests a natural way to construct a controller: by selecting, at each state $\bx$, the weighted centroid of the feasible control map $\Uc(\bx)$.

A weighted centroid $\map{\bmu}{\real^n}{\real^m}$ of a feasible control map~$\setmap{\Uc}{\real^n}{\real^m}$ is defined as:
\begin{equation}\label{eq:wc}
\bmu(\bx; \Uc) = \frac{\int_{\Uc(\bx)} \bu p(\bx,\bu) d\bu}{\int_{\Uc(\bx)} p(\bx,\bu) d\bu},
\end{equation}
where $\map{p}{\real^n \times \real^m}{\realpos}$ is a positive weight function that assigns a weight  to each control $\bu$ at state $\bx$. Our intuition is that if the boundaries of the set change smoothly, then the weighted centroid will change smoothly as well. While this idea may seem simple, formally developing the associated theory requires care, which we begin in this section by focusing on systems with a single input.

\begin{remark}\label{rmk:support}\longthmtitle{Support of Weight Functions}
For the weighted centroid~\eqref{eq:wc} to be well defined, the denominator must be nonzero.
Since we only place mild structural assumptions on~$\Uc$,
we restrict the weight function $p$ to be strictly positive for all $(\bx,\bu)$, i.e., its support is $\real^m$, for any $\bx$. This restriction can be relaxed if there exists a known region $\Wc\subset\real^m$ such that $\Wc\cap\Uc(\bx)$ is nonempty for all $\bx$, in which case the support (in $\bu$) of $p$ need only be $\Wc$. This is the case, for example, when we have CLFs and CBFs for systems with bounded inputs.~\hfill$\bullet$
\end{remark}

\subsection{Single-Input Case}

To begin our exposition, we consider the case when system \eqref{sys:ctrl-affine} has a single input ($m=1$).
 For this one-dimensional case, regardless of the number of constraints~$q$,
 the domain of integration is simply an interval on the real line:
\begin{equation}\label{eq:feas_single}
\Uc(\bx) \defeq \cap_{i\in\until{q}} \Omega_i(\bx)= \setdefb{u\in \real}{c(\bx)\leq u \leq d(\bx)}.
\end{equation}
for endpoint functions $\map{c,d}{\real^n}{\real}$. The following result characterizes the smoothness of the weighted centroid~\eqref{eq:wc}.

\begin{lemma}\longthmtitle{Single Input - Bounded}\label{lem:single_input}
Consider the weighted centroid~\eqref{eq:wc} of the admissible control map~\eqref{eq:feas_single} with a single input~$m=1$. Suppose that:
\begin{enumerate}
\item the functions $c$, $d$, and $p$ are  $\Cc^{\ell}$ at $\bx^*$;
\item the maps $\{\Omega_i\}_{i=1}^q$ are strictly compatible at $\bx^*$.
\end{enumerate}
Then $\mu$ is $\Cc^\ell$ at $\bx^*$.
\end{lemma}
\begin{proof}
First, we note that strict compatibility guarantees $c(\bx^*)\neq d(\bx^*)$ at $\bx^*$. Since $c$ and $d$ are continuous, there exists $\delta>0$ such that $c(\bx)\neq d(\bx)$ for all $\bx\in\Ball{\delta}{\bx^*}$. Thus, the function~$\mu$ is well-defined on $\Ball{\delta}{\bx^*}$. To establish that the partial derivatives of $\mu$ exist and are continuous on $\bx\in\Ball{\delta}{\bx^*}$, we apply the Leibniz integral rule \cite{HF:73} to show:
\begin{align}\label{eq:Leibniz}
\frac{\partial}{\partial x_i} \int_{c(\bx)}^{d(\bx)} \psi(\bx&,u) du = \int_{c(\bx)}^{d(\bx)} \frac{\partial \psi}{\partial x_i} (\bx,u) du\\&+\psi(\bx,d(\bx))\frac{\partial d}{\partial x_i}(\bx) -\psi(\bx,c(\bx))\frac{\partial c}{\partial x_i}(\bx). \nonumber
\end{align}
where $\psi$ is any smooth function. Because the functions $c$, $d$, and $\psi$ are smooth, each term on the right hand side exists and is continuous. Thus, using the Leibniz integral rule with $p$ and $u\cdot p$ as $\psi$, $\mu$ is continuously differentiable via the quotient rule of differentiation. Note that under the quotient rule, the denominator in each derivative remains simply some power of the original denominator, so it never takes the value of zero. We can then continue to use the Leibniz integral rule in order to show the existence and continuity of higher-order partial derivatives. This concludes the proof.
\end{proof}

Lemma \ref{lem:single_input} confirms our intuition on the smoothness of the weighted centroid selection function for the single-input case: $\mu$ inherits the smoothness properties from the endpoints $c$ and $d$. Although the proof using the Leibniz integral rule does not generalize to higher input dimensions, it reveals some issues that we shall discuss before proceeding to more general cases.

\subsection{Proper Weight Functions}
Lemma \ref{lem:single_input} does not apply when the integral bound $c(x)=-\infty$ or $d(x)=\infty$ becomes unbounded, which can happen when $b_i(x)=0$ for some $i\in\until{q}$. In fact, neither of the functions $c$ nor $d$ are well-defined there, let alone enjoy a notion of continuity or smoothness. In addition, when the interval of integration is unbounded, there is a question of whether the integrals in~\eqref{eq:wc}  converge. We will address these issues by carefully choosing the weight function~$p$.

Without loss of generality (through scaling), the weight functions $p$ that ensure the integrals defining the weighted centroid~\eqref{eq:wc} converge are probability density functions (PDFs). From this perspective, our weighted centroid is simply the expected value over the distribution, truncated by~\eqref{eq:ctrl-affine-feas}.

\begin{example}[Gaussian Formula]
The first natural choice to consider is the normal distribution PDF:
\begin{equation}
\label{eq:gaussian_one}
p^1_\gau(\bx,u) =  ((2\pi)^{1/2}\sigma(\bx))^{-1}\exp\Big(-\frac{1}{2}\frac{u^2}{\sigma(\bx)^2}\Big),
\end{equation}
where $\map{\sigma}{\real^n}{\real}$ represents the (potentially state-dependent) standard deviation. Here we consider the single-input case $m=1$. The expected value (in $u$) of a truncated Gaussian distribution is given by:
$$
\mu_\gau(\bx;\Uc) = \sigma(\bx)\frac{p^1_\gau(c(\bx)/\sigma(\bx))-p^1_\gau(d(\bx)/\sigma(\bx))}{\Phi(d(\bx)/\sigma(\bx))-\Phi(c(\bx)/\sigma(\bx))},
$$
where $\Phi$ is the standard cumulative distribution function (zero-mean and unit variance). This solution confirms the conclusion of Lemma~\ref{lem:single_input}. In addition, its derivatives are well-defined and converge smoothly as $c$ or $d$ go unbounded. ~\hfill $\bullet$
\end{example}

\subsection{Single Inequality Constraint}\label{sec:single_ineq}
We now consider the case of a single constraint ($q=1$) with multiple inputs ($m\geq 1$). Here, the admissible control set is defined by a single affine inequality ${a_1(\bx)+\bb_1(\bx)\bu \leq 0}$. A common approach, cf. \cite{EDS:89a}, is to restrict the input via a parametrization, $\bu=\bb_1^\top(\bx)v$, so that it lies in the direction of the constraint normal $\bb_1^\top(\bx)$. Effectively, the restriction reduces the control affine inequality to:
$$
a_1(\bx)+\bb_1(\bx)\bu \leq 0 \xRightarrow{\bu=\bb_1^\top(\bx)v}  a_1(\bx)+\|\bb_1(\bx)\|^2v \leq 0,
$$
which is now a scalar inequality in $v\in\real$.

\begin{example}[Sontag's Formula]\label{ex:sontag}
To apply the weighted centroid methodology in this reduced problem, we choose a weight function for the new scalar variable $v$:
    $$
    p_\son(\bx,v) = (1/2)(v^2+1)^{-3/2}.
    $$
We may evaluate the
weighted centroid~\eqref{eq:wc} as:
\begin{align*}
\mu_{\son,v}(\bx;\Omega_1)&=\frac{\int_{-\infty}^{d(\bx)} v p_\son(\bx,v)dv}{\int_{-\infty}^{d(\bx)}  p_\son(\bx,v)dv}
= \frac{-\frac{1}{2\sqrt{d(\bx)^2+1}}}{\frac{1}{2}\big(\frac{d(\bx)}{\sqrt{d(\bx)^2+1}}+1\big)}\\
 &= \frac{-1}{d(\bx)+\sqrt{d(\bx)^2+1}}= d(\bx)-\sqrt{d(\bx)^2+1}.
\end{align*}
where $d(\bx)=-a_1(\bx)/\|\bb_1(\bx)\|^2$.
We then convert this to the original input coordinate with $\bu=\bb_1^\top(\bx)v$ to obtain:
\begin{equation}\label{eq:sontag}
    \bmu_\son(\bx;\Omega_1) = \tfrac{-a_1(\bx)-\sqrt{a_1(\bx)^2+\|\bb_1(\bx)\|^4}}{\|\bb_1(\bx)\|^2}\bb_1^\top(\bx),
\end{equation}
which is Sontag's universal formula for stabilization, cf.~\cite{EDS:89a}, that constructively provides an analytic controller for Artstein's theorem, albeit limited to the one-constraint case. This derivation shows that Sontag's formula can be interpreted as a weighted centroid of the admissible control set.~\hfill $\bullet$
\end{example}
While Sontag's formula is known to be smooth, this property follows from its original derivation using tools such as the Implicit Function Theorem (see also~\cite{MHC-PO-GB-ADA:23-csl}). In contrast, within the weighted centroid framework, smoothness is not immediate due to the possibility of state-dependent integration endpoint $d(\bx)$ going unbounded, e.g., as $\bb_1(\bx)\rightarrow \bzero$. In particular, establishing smoothness requires understanding how differentiation interacts with the integrals. To this end, we introduce the notion of \emph{integrable differentiability}.

\begin{definition}
   Let $\map{\psi}{\real^n\times\real}{\real}$ be a $\Cc^\ell$ function. Denote by $\Dc^\ell(\psi)$
   the collection of all of its partial derivatives
   up to order $\ell$, including $\psi$ itself. We say that $\psi$ is \emph{$\ell$-integrably differentiable}
   if, for any $\phi\in\Dc^\ell$,
   the integral:
   \begin{align}\label{eq:integrable_diff}
\int_{-\infty}^\infty |\phi(\bx,v)|dv
\end{align}
    converges for each $\bx$.~\hfill$\diamond$
\end{definition}
Integrable differentiability enables the following result.

\begin{theorem}\label{thm:Leibniz_converging}
\longthmtitle{Single Input - Unbounded}
Consider a single (${q=1}$) admissible control map~\eqref{eq:ctrl-affine-feas}. Suppose~that:
\begin{enumerate}
    \item the functions $a_1$, $\bb_1$, and $p$ are $\Cc^\ell$ on $\Sc_1$;
    \item the map $\Omega_1$ is strictly compatible on $\Sc_1$;
    \item both functions $(\bx,v)\mapsto p(\bx,v)$ and $(\bx,v)\mapsto vp(\bx,v)$ are $\ell$-integrably differentiable in $\bx$;
\item the following limit converges for all $\phi\in \Dc^{\ell-1}(p)$:
\begin{equation}\label{eq:extreme-limit}
\lim_{v\rightarrow \infty} v^{\ell+2}\phi(\bx,v)
\end{equation}
\end{enumerate}
Then, the controller:
\begin{equation}\label{eq:single-input-controller-and-map}
\begin{aligned}
    \bk(\bx) &=  \bb_1(\bx)\T \mu(\bx;\Omega_{v}), \\
    \Omega_{v}(\bx) &=  \setdefb{v\in\real}{-\infty < v \leq -\tfrac{a_1(\bx)}{\|\bb_1(\bx)\|^2}},
\end{aligned}
\end{equation}
with $\mu$ the weighted centroid function given by~\eqref{eq:wc},
is a $\Cc^\ell$ function on $\Sc_1$ and satisfies:
$$
a_1(\bx)+\bb_1(\bx)\bk(\bx)\leq 0,\quad \forall \bx\in \Sc_1.
$$
\end{theorem}

\begin{proof}
    Our proof focuses on first-order differentiability with $\ell=1$. This same process can be repeated for higher-order derivatives using the assumptions outlined in the  statement. Following similar steps to the proof of Lemma \ref{lem:single_input}:
    \begin{equation}\label{eq:partial-integral-psi}
    \begin{aligned}
        \pdv{}{x_i}\int_{-\infty}^{d(\bx)}\psi(\bx,v)dv = & \int_{-\infty}^{d(\bx)}\pdv{\psi}{x_i}(\bx,v)dv \\ & + \psi(\bx,d(\bx))\pdv{d}{x_i}(\bx),
    \end{aligned}
    \end{equation}
    where $d(\bx)\coloneqq -\tfrac{a_1(\bx)}{\beta(\bx)}$ and $\beta(\bx)\coloneqq \|\bb_1(\bx)\|^2$. We now illustrate that the above expression is continuous for all $\bx\in\Sc_1$. We consider two separate cases.

    \underline{$\bb_1(\bx)\neq\bzero$:} In this case, $d$ is continuous at $\bx$, and is locally
    lower bounded by some constant $c\in\real$.
    The integral over the domain $(-\infty, d(\bx)]$ is the sum of the integral over the domains $(-\infty,c]$ and $[c,d(\bx)]$. From these summands, the former is a finite constant from integrable differentiability. As for the latter, the integral adopts continuity from all the functions involved. Therefore, the sum is continuous.

    \underline{$\bb_1(\bx)=\bzero$:} We establish differentiability of $\bk$ by considering the limiting behavior as $\bb_1(\bx)\rightarrow\bzero$ or, equivalently, as $d(\bx)\rightarrow\infty$ since we must have $a_1(\bx) < 0$ as $\bb_1(\bx)\rightarrow \bzero$ from the strict compatibility assumption.

    First, by integrable differentiability, for $\psi(\bx,v)=p(\bx,v)$ and $\psi(\bx,v)=vp(\bx,v)$, the term:
    \begin{equation*}
        \int_{-\infty}^{\infty}\pdv{\psi}{x_i}(\bx,v)dv = \lim_{d(\bx)\rightarrow\infty}\int_{-\infty}^{d(\bx)}\pdv{\psi}{x_i}(\bx,v)dv,
    \end{equation*}
    is well defined. Next, note that with $d(\bx)=-\tfrac{a_1(\bx)}{\beta(\bx)}$ we have:
    \begin{equation*}
        \begin{aligned}
            \pdv{d}{x_i}(\bx) =  \frac{-1}{\beta(\bx)}\pdv{a_1}{x_i}(\bx)  + \pdv{\beta}{x_i}(\bx)\frac{a_1(\bx)}{\beta^2(\bx)}.
        \end{aligned}
    \end{equation*}
    Now, let $\bx_0$ be any point satisfying $\bb_1(\bx_0)=\bzero$ and note that since $a_1(\bx)<0$ whenever $\bb_1(\bx)=\bzero$, there exists a $\varepsilon>0$ such that $a_1(\bx)<0$ for all $\bx\in\mathcal{B}_{\varepsilon}(\bx_0)$. Hence, for $\bx\in\mathcal{B}_{\varepsilon}(\bx_0)\setminus\{\bx_0\}$, we have:
    \begin{equation*}
        \begin{aligned}
            \psi(\bx,d(\bx))&\pdv{d}{x_i}(\bx) \\
            =& \frac{\psi(\bx,d(\bx))}{a_1(\bx)}\bigg[\frac{-a_1(\bx)}{\beta(\bx)}\pdv{a_1}{x_i}(\bx)  + \frac{a_1^2(\bx)}{\beta^2(\bx)}\pdv{\beta}{x_i}(\bx)\bigg] \\
            =& \frac{\psi(\bx,d(\bx))}{a_1(\bx)}\bigg[d(\bx)\pdv{a_1}{x_i}(\bx)  + d^2(\bx)\pdv{\beta}{x_i}(\bx)\bigg].
        \end{aligned}
    \end{equation*}
    Then since the limits~\eqref{eq:extreme-limit} exists, the following limit:
    \begin{equation*}
        \begin{aligned}
            \lim_{\bx\rightarrow\bx_0}\psi(\bx,d(\bx))d^j(\bx),\quad j\in\{1,2\},
        \end{aligned}
    \end{equation*}
    with $\psi(\bx,v)=p(\bx,v)$ or $\psi(\bx,v)=vp(\bx,v)$ also exists. Consequently, with $a_1$ and $\bb_1$ being continuous at $\bx$, each term on the right hand side in~\eqref{eq:partial-integral-psi} converges as $\bx\rightarrow \bx_0$. That is, all (first-order) partial derivatives are well-defined as $\bb(\bx)\rightarrow\bzero$. Since all partial derivatives exist and are continuous on any open set such that $\bb_1(\bx)\neq\bzero$ (from Case (i)) and the limit of these partial derivatives exists as $\bx$ tends to points such that  $\bb_1(\bx)=\bzero$, it follows that all partial derivatives exists and are continuous at such points as well, thereby ensuring that $\mu$ as defined in \eqref{eq:wc} is continuously differentiable. For higher-order derivatives, we can take the derivative of \eqref{eq:partial-integral-psi} using the product rule and the Leibniz rule of differentiation, and apply the set of assumptions outlined in the statement to establish the existence of such higher-order derivatives.

    Since $\mu$ is a weighted centroid of the convex set $\Omega_{v}(\bx)$, it satisfies $\mu(\bx;\Omega_v)\in\Omega_{v}(\bx)$ for all $\bx\in\Sc_1$ \cite[p. 25]{SB-LV:09}.
    Hence, for all $\bx\in\Sc_1$, we can show the controller $\bk$ satisfies:
    \begin{equation*}
        \begin{aligned}
             a_1(\bx) + \bb_1(\bx)\bk(\bx) =  a_1(\bx) + \mu(\bx,\Omega_{v})\|\bb_1(\bx)\|^2.
        \end{aligned}
    \end{equation*}
    When $\bb_1(\bx)=\bzero$ we have $a_1(\bx) < 0$ by assumption, and when $\bb_1(\bx)\neq\bzero$ we have:
    \begin{equation*}
         a_1(\bx) + \bb_1(\bx)\bk(\bx) \leq a_1(\bx) - \frac{a_1(\bx)}{\|\bb(\bx)\|^2}\|\bb_1(\bx)\|^2 = 0,
    \end{equation*}
    for all $\bx\in\Sc_1$, as desired.
\end{proof}

Theorem~\ref{thm:Leibniz_converging} outlines assumptions under which each term in the Leibniz integral rule \eqref{eq:partial-integral-psi} converges. For differentiability with $\ell=1$, the assumptions are not too
restrictive.
In addition to regularity assumptions,
~\eqref{eq:extreme-limit} imposes that the original probability distribution is not too tail-heavy. One may verify that these conditions hold for
~\eqref{eq:sontag-like_weight} and the Gaussian distribution~\eqref{eq:gaussian_one}. Beyond first derivatives, the number of terms whose convergence must be verified grows combinatorially, making such conditions increasingly difficult to verify in practice.

Note, however, that the assumptions made in Theorem \ref{thm:Leibniz_converging} are not necessary conditions. For instance, one may verify that, beyond first derivatives, individual terms derived from the Sontag's weight function~$p_\son$ do not converge, yet the weighted centroid associated with Sontag's formula is known to be an analytic function~\cite{EDS:89a}. On the other hand, weight functions with sufficiently fast decay, e.g., from the exponential family, ensure convergence of all terms due to the dominance of exponential tails over polynomials.

\begin{example}[Sontag-like Formulas]\label{ex:sontag-like}
Theorem~\ref{thm:Leibniz_converging} provides insight into choosing weight functions:
condition~\eqref{eq:extreme-limit} requires sufficient decay of the weight function at infinity.
Motivated by this, we consider a family of weight functions of the form:
\begin{equation}\label{eq:sontag-like_weight}
    p(\bx,v) = (v^2+1)^{-r/2},
\end{equation}
with $r\geq 3$.
For example, with $r=5$ we have:
\begin{equation*}
    \bmu(\bx;\Omega_1) = \frac{-\bb_1^\top(\bx)}{d(\bx)(2d(\bx)^2+3)+2(d(\bx)^2+1)^{3/2}},
\end{equation*}
and with $r=7$ the weighted centroid becomes:
\begin{equation*}
    \bmu(\bx;\Omega_1) = \frac{-3\bb^\top_1(\bx)}{8d(\bx)^5+20d(\bx)^3+8(d(\bx)^2+1)^{5/2}+15d(\bx)}.
\end{equation*}
As $r$ increases,
\eqref{eq:extreme-limit} is satisfied for higher values of $\ell$, improving the order of differentiability of the resulting formulas. As $r\to\infty$,
\eqref{eq:sontag-like_weight}
approaches the Gaussian weight, which satisfies the assumptions for all orders due to its exponential decay.~\hfill $\bullet$
\end{example}

\subsection{Mirrored Constraints}
The preceding result allows for constructing smooth controllers satisfying a single affine inequality constraint. While controllers that can achieve this have been around for quite a while (e.g., dating back to Sontag's original formulation in \cite{EDS:89a}), the weighted centroid perspective taken herein allows for extending such ideas to the multi-constraint case.  In particular, without any modification, the results of this section may be directly used to construct smooth controllers satisfying a class of constraints we refer to as \emph{mirrored}.
\begin{definition}\longthmtitle{Mirrored Constraints}
    A pair of affine constraints $\{a_i(\bx) + \bb_i(\bx)\bu\leq0\}_{i\in\until 2}$ is \emph{mirrored} for \eqref{sys:ctrl-affine} if the constraints normals satisfy $\bb_2(\bx)= -\bb_1(\bx)$ for all $\bx\in\R^n$.~\hfill $\diamond$
\end{definition}
When a pair of constraints is mirrored, we may parameterize the input as $\bu = \bb_1(\bx)\T v$, $v\in\R$, to recast the constraints as:
\begin{equation*}
\begin{aligned}
    a_1(\bx) + \|\bb_1(\bx)\|^2v \leq & 0 \\
    a_2(\bx) - \|\bb_1(\bx)\|^2 v \leq & 0,
\end{aligned}
\end{equation*}
leading to the admissible control map:
\begin{equation}\label{eq:mirrored-control-map}
    \Omega_v(\bx) = \bigg\{v\in\R\,:\underbrace{\frac{a_2(\bx)}{\|\bb_1(\bx)\|^2}}_{c(\bx)} \leq v \leq \underbrace{\frac{-a_1(\bx)}{\|\bb_1(\bx)\|^2}}_{d(\bx)}\bigg\},
\end{equation}
whose weighted centroid yields a smooth controller satisfying the original pair of constraints.

\begin{proposition}\longthmtitle{Single Mirrored Pair}\label{prop:mirrored-constraints}
    Consider a pair of affine constraints $a_i(\bx) + \bb_i(\bx)\bu\leq0$, $i\in\{1,2\}$. that are mirrored and strictly compatible on $\mathcal{S}=\mathcal{S}_1\cap\mathcal{S}_2$. Provided $a_i$, $\bb_i$, $i\in\{1,2\}$ are $\mathcal{C}^\ell$ and $\bb_i(\bx)\neq\bzero$ on $\mathcal{S}$, then the controller $\bk(\bx)=\mu(\bx;\Omega_{v})\bb_1(\bx)^\top$, with $\mu$ the weighted centroid function given by~\eqref{eq:wc} of the set $\Omega_v$ as in \eqref{eq:mirrored-control-map}, is $\mathcal{C}^\ell$ on $\mathcal{S}$ and satisfies $a_i(\bx) + \bb_i(\bx)\bk(\bx)\leq0$ for $i\in\{1,2\}$.
\end{proposition}
\begin{proof}
    Since $a_i$ and $\bb_i$ are $\mathcal{C}^\ell$ and $\bb_i(\bx)\neq0$, it follows that $c$ and $d$ as defined in \eqref{eq:mirrored-control-map} are $\mathcal{C}^\ell$. The claim that resulting controller is $\mathcal{C}^\ell$ then follows directly from Lemma~\ref{lem:single_input}. Furthermore, since $\mu(\bx,\Omega_{v})$ is the weighted centroid of a convex set, we have $\mu(\bx,\Omega_{v})\in\Omega_{v}$ and constraint satisfaction follows a similar argument to that in Theorem \ref{thm:Leibniz_converging}.
\end{proof}

\begin{proposition}\label{prop:mirrored-constraints-b=0}
\longthmtitle{Single Mirrored Pair with vanishing normals}
    Let the conditions of Proposition \ref{prop:mirrored-constraints} hold with the exception of the assumption that $\bb_i(\bx)\neq\bzero$ for all $\bx\in\mathcal{S}$ (i.e., there may exist points where $\bb_i(\bx)=\bzero$).
    Suppose the assumptions of Theorem~\ref{thm:Leibniz_converging} hold with the changes of strict compatibility to be between the maps $\{\Omega_1, \Omega_2\}$ on $\Sc$ and the limit~\eqref{eq:extreme-limit} to also hold as $v\rightarrow -\infty$.
    Then, the controller $\bk(\bx)=\mu(\bx;\Omega_{v})\bb_1(\bx)^\top$, with $\Omega_v$ as in \eqref{eq:mirrored-control-map}, is $\mathcal{C}^1$ on $\mathcal{S}$ and satisfies $a_i(\bx) + \bb_i(\bx)\bk(\bx)\leq0$ for $i\in\{1,2\}$.
\end{proposition}
\begin{proof}
    The proof follows the same steps as those in Theorem~\ref{thm:Leibniz_converging} and is omitted for brevity.
\end{proof}

Proposition~\ref{prop:mirrored-constraints-b=0} is our first result towards extending our single constraint results to the multi-constraint setting.
We now extend Sontag's formula to
mirrored constraints.

\begin{example}[Mirrored Sontag's formula]
Using the results of this section, Sontag's formula from Example \ref{ex:sontag} may be extended to handle mirrored constraints.
We have:
\begin{equation*}
    \int_{c(\bx)}^{d(\bx)} p_\son(\bx,v)dv = \frac{1}{2}\left(\frac{d(\bx)}{\sqrt{d(\bx)^2+1}}-\frac{c(\bx)}{\sqrt{c(\bx)^2+1}}\right),
\end{equation*}
\begin{equation*}
    \int_{c(\bx)}^{d(\bx)} vp_\son(\bx,v)dv = \frac{1}{2}\left(\frac{1}{\sqrt{c(\bx)^2+1}}-\frac{1}{\sqrt{d(\bx)^2+1}}\right),
\end{equation*}
yielding the weighted centroid for the case $\bb_1(\bx)\neq \bzero$ and $\bb_2(\bx)\neq \bzero$, in the $v$-coordinate as:
\begin{equation}
    \mu_{\mir,v}(\bx;\Omega_v) =
        \frac{\frac{1}{\sqrt{c(\bx)^2+1}}-\frac{1}{\sqrt{d(\bx)^2+1}}}{\frac{d(\bx)}{\sqrt{d(\bx)^2+1}}-\frac{c(\bx)}{\sqrt{c(\bx)^2+1}}},
\end{equation}
which is
continuously differentiable by Proposition~\ref{prop:mirrored-constraints-b=0} with $\ell=1$. We now derive the formula in $\bu$-coordinates as:
\begin{equation}\label{eq:sontag_mirror}
    \bmu_\mir(\bx;\Uc) = \begin{cases}
        \bzero, &~\bb_1(\bx)= \bb_2(\bx) = \bzero \\
        \bmu_\son(\bx;\Omega_1), &~\bb_2(\bx) = \bzero \\
        \bmu_\son(\bx;\Omega_2), &~\bb_1(\bx) = \bzero \\
        \mu_{\mir,v}(\bx;\Omega_v)\bb_1(\bx)^\top, &~\text{otherwise}
    \end{cases}
\end{equation}
where $\bmu_\son$ is defined in \eqref{eq:sontag}.
~\hfill$\bullet$
\end{example}
Proposition~\ref{prop:mirrored-constraints-b=0} is limited to only two constraints that are mirrored due to an obstruction we discuss next.

\subsection{Nonsmoothness from Constraint Switching}
The smoothness requirement of the endpoint functions $c$ and $d$ in Lemma~\ref{lem:single_input} restricts the scope of our result. Even when all biases $a_i$'s and normals $\bb_i$'s are smooth, $c$ and $d$ may not be smooth. For example, consider the following collection of admissible control maps with ${\bb_i(\bx) \equiv 1}$:
\begin{align*}
\Omega_1(\bx) &= \setdefb{u\in\real}{u\leq -a_1(\bx)}\\
\Omega_2(\bx) &= \setdefb{u\in\real}{u\leq -a_2(\bx)}\\
\implies \Uc(\bx) &= \setdefb{u\in\real}{u\leq \min\{-a_1(\bx),-a_2(\bx)\}}
\end{align*}
Here, $d(\bx)=\min\{-a_1(\bx),-a_2(\bx)\}$, which is not necessarily smooth if $a_1(\bx)=a_2(\bx)$. This example illustrates how nonsmoothness can arise from switching in constraints that define the boundary of the feasible set. Such switching can occur when some of the constraint normals are linearly dependent.

The loss of smoothness is fundamentally driven by a constraint regularity issue, rather than simply the smoothness of individual constraint bias and normal functions. Such questions are classically studied in optimization via constraint qualification conditions. This connection is especially relevant in control, where optimization-based controllers are used to satisfy CLF and CBF constraints. In such settings, \textit{linear independence constraint qualification (LICQ)}~\cite{JN-SW:06}
guarantees the local Lipschitzness of the controller~\cite{PM-AA-JC:25}. Although the weighted centroid construction does not rely on solving an optimization problem, the same geometric obstruction involving linear dependency appears,
motivating our
focus
in the next section on constraints with linearly independent normals.

\section{General Multiple-Constraint Cases}\label{sec:general-case}
In this section, we generalize the weighted centroid methodology to an arbitrary number of admissible control maps. The main difficulty in establishing smoothness of the weighted centroid comes from the dependency of the domain of integration $\Uc(\bx) = \cap_{i\in\until{q}} \Omega_i(\bx)$  on the state~$\bx$. In the single-input case, this dependence is handled using the Leibniz integral rule to analyze the derivatives, but such an approach does not extend easily to higher-dimensional input spaces.

We begin with the case where the constraint normals $\bb_i$ are linearly independent, allowing us to effectively decouple the constraints through a local coordinate transformation and apply the previous results. Then, we turn to the subtler linearly dependent case, where degeneracies arise, and show that we can recover smoothness with Gaussian weight functions.

\subsection{Linearly Independent Constraint Normals}

For a state $\bx^*$, let $\Ac$ be a set of indices where the constraint normals $\{\bb_i(\bx^*)\}_{i\in\Ac}$ are linearly independent. Without loss of generality, relabel the indices so
$\Ac=\{1,\dots,|\Ac|\}$. Define:
$$
\ba(\bx^*) = {\tiny\begin{bmatrix}a_1(\bx^*) \\\vdots \\ a_{|\Ac|}(\bx^*)\end{bmatrix}},\quad\bB(\bx^*) = {\tiny\begin{bmatrix}\bb_1(\bx^*) \\\vdots \\ \bb_{|\Ac|}(\bx^*)\end{bmatrix}}.
$$
Since the constraint normals $\{\bb_i(\bx^*)\}_{i\in\Ac}$ are linearly independent, we can find a subset of standard unit vectors $\{\hat \be_j\}_{|\Ac|+1\leq j \leq m}$ (relabeled without loss of generality) that completes the basis of $\real^m$ at $\bx^*$. By continuity, the row vectors $\{\bb_i(\bx^*)\}_{i\in\Ac}$ together with $\{\hat \be_j\}_{|\Ac|+1\leq j \leq m}$ remain linearly independent locally on a sufficiently small ball $\Ball{\delta}{\bx^*}$.

Applying the Gram-Schmidt orthogonalization process to the vectors $\{\hat \be\}_{|\Ac|+1\leq i \leq m}$ relative to the constraint normals $\{\bb_i(\bx)\}_{i\in\Ac}$ for $\bx\in \Ball{\delta}{\bx^*}$, we obtain smooth orthogonal vectors $\{\bT_i(\bx;\bx^*)\}_{|\Ac|+1\leq i \leq m}$. This yields a smooth local coordinate transformation matrix:
$$
\bT(\bx;\bx^*) = \small\begin{bmatrix} \bB(\bx) \\ \bT_{|\Ac|+1}(\bx;\bx^*) \\ \vdots \\ \bT_m(\bx;\bx^*)\end{bmatrix},
$$
for
the following input transformation into a new coordinate:
\begin{equation}\label{eq:w-coordinates}
\bw  =  {\small \begin{bmatrix}\ba(\bx) \\ \bzero \end{bmatrix}}+\bT(\bx;\bx^*)\bu.
\end{equation}
Although the transformation is only valid locally, it is sufficient for establishing smoothness of the weighted centroid at $\bx$.
Particularly, the constraint map $\Uc_\Ac(\bx)=\cap_{i\in\Ac}\Omega_i(\bx)$ is constant in the state~$\bx$ when expressed in the new coordinate:
$$
\Omega_0 = \setdefb{\bw\in\real^m}{w_i\leq 0,~\forall i\in\Ac} ,
$$
for each $\bx\in \Ball{\delta}{\bx^*}$. Hence, we may express an integral of any smooth function $\map{\psi}{\real^n \times \real^m}{\real}$
in $\bw$ coordinates as
\begin{multline*}
\int_{\Uc_\Ac(\bx)}\psi(\bx,\bu)d\bu=\\
\int_{\Omega_0} \psi\left(\bx,\bT^{-1}(\bx;\bx^*)\Big(\bw- \begin{bmatrix}\ba(\bx) \\ \bzero \end{bmatrix}\Big) \det(\bT^{-1}(\bx;\bx^*))\right)\!d\bw.
\end{multline*}
With a constant domain of integration, the analysis of the smoothness of the weighted centroid simplifies as we can rely on the dominated convergence theorem to swap the order of differentiation and integration. As a result, we derive the following sufficient conditions that ensures smoothness.

\begin{proposition}\longthmtitle{Smoothness from Linear Independence}\label{prop:lin_ind}
Let $\Ac\subseteq\until{q}$ be a subset of indices of the admissible control maps $\{\Omega_i\}_{i\in\until{q}}$
in~\eqref{eq:ctrl-affine-feas}.
Consider the weighted centroid~\eqref{eq:wc} of the intersection $\Uc_\Ac(\bx)=\cap_{i\in\Ac}\Omega_i(\bx)$ of admissible control maps with only indices from
$\Ac$. Let $\bx^* \in \real^n$.
Assume that:
\begin{enumerate}
\item the functions $\{a_i\}_{i=1}^\Ac$, $\{\bb_i\}_{i=1}^\Ac$ and $p$, are $\Cc^{\ell}$;
\item for every partial derivative $\phi \in \Dc^\ell(p)$ of $p$,
the integral:
\begin{align}\label{eq:infinite-integral-converges}
\int_{\real^m} \Pc^{\ell+1}(\bu)\phi(\bx,\bu)d\bu
\end{align}
converges for all $\bx$, where $\Pc^{\ell+1}(\bu)$ denotes any $\ell+1$-th order polynomial of $\bu$;
\item the constraint normals $\{\bb_i(\bx^*)\}_{i\in\Ac}$ are linearly independent.
\end{enumerate}
Then $\bmu(\cdot;\Uc_\Ac)$ is $\Cc^\ell$ at $\bx^*$. Thus, if the assumptions are satisfied with $\Ac=\until{q}$, the weighted centroid $\bmu(\cdot;\Uc)$ is $\Cc^\ell$ at~$\bx^*$.
\end{proposition}
\begin{proof}
To compute the derivatives of $\bmu$ defined in \eqref{eq:wc}, we express the integrals in the $\bw$ coordinates where we may apply Lebesgue's dominated convergence theorem \cite[Thm. 1.34]{WR:87} to differentiate under the integral sign, after which we can convert back to the $\bu$ coordinates as follows:
    \begin{align}
&\frac{\partial}{\partial x_i} \int_{\Uc_\Ac(\bx)}\psi(\bx,\bu)d\bu \nonumber\\
 &=\frac{\partial}{\partial x_i}\left(\det(\bT^{-1}(\bx))\right) \int_{\Omega_0} \psi\Big(\bx,\bT^{-1}(\bx)\big(\bw- \begin{bmatrix}\ba(\bx) \\ \bzero \end{bmatrix}\big)\Big)d\bw \nonumber\\
  &~+ \det(\bT^{-1}(\bx)) \Big[\int_{\Omega_0}\frac{\partial\psi}{\partial x_i}\!\Big(\bx,\bT^{-1}(\bx)\big(\bw- \begin{bmatrix}\ba(\bx) \\ \bzero \end{bmatrix}\big)\Big)d\bw \nonumber\\
  &~+\int_{\Omega_0}\frac{\partial\psi}{\partial \bu}\Big(\bx,\bT^{-1}(\bx)\big(\bw- \begin{bmatrix}\ba(\bx) \\ \bzero \end{bmatrix}\big)\Big)\nonumber\\\
  &~\times\left(\frac{\partial \bT^{-1}}{\partial x_i}(\bx)\big(\bw- \begin{bmatrix}\ba(\bx) \\ \bzero \end{bmatrix}\big)-\bT^{-1}(\bx) \begin{bmatrix}\frac{\partial \ba}{\partial x_i}(\bx) \\ \bzero \end{bmatrix}\right)d\bw\Big] \nonumber\\
&= \frac{1}{\det(\bT^{-1}(\bx))} \frac{\partial}{\partial x_i}\left(\det(\bT^{-1}(\bx))\right) \int_{\Uc_\Ac(\bx)} \psi(\bx,\bu)d\bu\nonumber\\
&~+ \int_{\Uc_\Ac(\bx)}\frac{\partial\psi}{\partial x_i}(\bx,\bu) d\bu+  \int_{\Uc_\Ac(\bx)}\frac{\partial\psi}{\partial \bu}(\bx,\bu)\nonumber\\
&~\times\left(\frac{\partial \bT^{-1}}{\partial x_i}(\bx)\bT(\bx)\bu-\bT^{-1}(\bx) \begin{bmatrix}\frac{\partial \ba}{\partial x_i}(\bx) \\ \bzero \end{bmatrix}\right)d\bu.
  \label{eq:integral_partials}
\end{align}
Using the identity above, we find that repeated partial derivatives in $\bx$ result in a sum of integrals of partial derivatives $\phi\in\Dc^\ell$ multiplying with polynomials in $\bu$ (whose coefficients depend smoothly on $\bx$). By the second assumption in the statement, these integrals converge to finite values. Hence, the partial derivatives exist and are continuous on $\Ball{\delta}{\bx^*}$.  Consequently, $\bmu$ is smooth at $\bx^*$.
\end{proof}

Given a weight function $p$ with sufficient regularity (e.g., a Gaussian weight), Proposition~\ref{prop:lin_ind} ensures the weighted centroid function~$\bmu$ is smooth at any state $\bx$ for which $\{\bb_i(\bx)\}_{i\in\until{q}}$ are linearly independent.
It is worth noting that this same assumption also makes it possible to select the ``vertex" control obtained by solving $\ba(\bx) + \bB(\bx)\bu=\bzero$, which is itself a smooth selection under linear independence. However, such a vertex control often requires unnecessarily large inputs or cancels beneficial dynamics when some $a_i(\bx)$ are already strictly negative. In contrast, the weighted centroid incorporates the distributional information encoded in $p$, biasing the controller toward more desirable regions of the admissible control set while still retaining smoothness.

Furthermore, the weighted centroid framework is flexible enough to handle certain linearly dependent cases where smoothness can still be recovered through appropriate coordinate transformation or weight-function design. We examine these structured cases next before addressing the fully degenerate situation in which some constraint normals $\bb_i(\bx)$ vanish.

\subsection{Structured Linearly Dependent Inequalities}
A common source of linear dependence arises when inequality constraints appear in \textit{mirrored} pairs, e.g., when certain quantities are constrained from both sides. Let $\Ac\subset \until{q}$ denote the index set of ``primary" constraints. We then assume that there is a subset $\Ic^+\subseteq \Ac$ such that each $i^+\in\Ic^+$ has a unique mirrored counterpart indexed by $i^- \notin \Ac$ such that:
$$
\bb_{i^-}(\bx) = -\bb_{i^+}(\bx),~\forall \bx\in \Sc.
$$
We collect all such mirrored indices in the set $\Ic^-$. Although each mirrored pair of constraints indexed $\{i^+,i^-\}$ are linearly dependent, the associated halfspaces $\{\Omega_{i^+}, \Omega_{i^-}\}$ have a useful geometric decomposition that allows us to avoid dealing with their intersection directly. This is formalized next.

\begin{lemma}\longthmtitle{Decomposition of Mirrored Halfspaces}
    \label{lem:mirrored_halfspaces}
    Consider a pair of affine constraints $a_i(\bx)+\bb_i(\bx)\bu \leq 0$, $i\in\{i^+,i^-\}$ and suppose that these constraints are mirrored and strictly compatible at $\bx$. Then the following identity holds:
    \begin{equation}\label{eq:decomposition_identity}
        \Omega_{i^+}(\bx)\cap\Omega_{i^-}(\bx) = \Omega_{i^+}(\bx)\setminus\Omega_{i^-}^c(\bx).
    \end{equation}
\end{lemma}
\begin{proof}
Strict compatibility of the mirrored pair implies:
$$
a_{i^-}(\bx) < \bb_{i^+}(\bx)\bu <  -a_{i^+}(\bx).
$$
Thus, the intersection of the complement sets $\Omega^c_{i^+}  \cap \Omega^c_{i^-}$ must be empty since, otherwise, there would exist a $\bu$ such that:
$$
a_{i^-}(\bx) \geq \bb_{i^+}(\bx)\bu \geq  -a_{i^+}(\bx),
$$
which would be a contradiction. Consequently, we have:
\begin{align*}
\Omega_{i^-}^c(\bx)&=  (\Omega_{i^+}(\bx)  \cap \Omega^c_{i^-}(\bx))\cup (\Omega^c_{i^+}(\bx)  \cap \Omega^c_{i^-}(\bx))\\
&= \Omega_{i^+}(\bx)  \cap \Omega^c_{i^-}(\bx),
\end{align*}
and therefore, we can decompose:
\begin{align*}
\Omega_{i^+}(\bx)&= (\Omega_{i^+}(\bx)  \cap \Omega_{i^-}(\bx))\cup (\Omega_{i^+}(\bx)  \cap \Omega^c_{i^-}(\bx))\\
&= (\Omega_{i^+}(\bx)  \cap \Omega_{i^-}(\bx))\cup  \Omega^c_{i^-}(\bx).
\end{align*}
Since $\Omega_{i^+}(\bx) \cap \Omega_{i^-}(\bx)$ is a subset of $\Omega_{i^-}(\bx)$, it is disjoint from $\Omega^c_{i^-}(\bx)$.
As such, we have the identity in \eqref{eq:decomposition_identity}.
\end{proof}

Lemma~\ref{lem:mirrored_halfspaces} allows us to exploit the mirrored structure to simplify the region of integration. The identity in~\eqref{eq:decomposition_identity} shows that the intersection of a mirrored pair can be expressed as the difference of two individual halfspaces. Consequently, an integral over the feasible control map $\Uc(\bx) = \bigcap_{i\in {\Ac}\cup{\Ic^-}}\Omega_i(\bx)$ can be equivalently decomposed into a combination (with additions and subtractions) of integrals, each taken over a region defined by only one constraint from each mirrored pair. As a result, we obtain the following smoothness property for weighted centroids of independent mirrored constraints.

\begin{theorem}
\longthmtitle{Independent Mirrored Pairs}
\label{thm:ind_mirror}
Let $\Ac\subseteq\until{q}$ be a subset of indices of the admissible control maps $\{\Omega_i\}_{i\in\until{q}}$
in~\eqref{eq:ctrl-affine-feas}, let $\Ic^+\subseteq\Ac$ be a subset of indices that have unique mirrored counterparts, and let the set $\Ic^-\subset \until{q}$ denote those mirrored indices. Consider the weighted centroid
\eqref{eq:wc} of the intersection $\Uc_{\Ac\cup\Ic^-}(\bx)=\cap_{i\in\Ac\cup\Ic^-}\Omega_i(\bx)$ of admissible control maps with only indices from
$\Ac\cup\Ic^-$.
Assume that:
\begin{enumerate}
\item the functions $a_i,\bb_i$, $\forall i\in\Ac\cup\Ic^-$, and $p$, are $\Cc^{\ell}$;
\item for every partial derivative $\phi \in \Dc^\ell(p)$ of $p$,
the integral:
\begin{align}\label{eq:infinite-integral-converges}
\int_{\real^m} \Pc^{\ell+1}(\bu)\phi(\bx,\bu)d\bu
\end{align}
converges for all $\bx$, where $\Pc^{\ell+1}(\bu)$ denotes any $\ell+1$-th order polynomial of $\bu$;
\item the vectors $\{\bb_i(\bx^*)\}_{i\in\Ac}$ are linearly independent.
\end{enumerate}
Then $\bmu(\cdot;\Uc_{\Ac\cup\Ic^-})$ is $\Cc^\ell$ at $\bx^*$. Consequently, if the assumptions are satisfied with $\Ac\cup\Ic^-=\until{q}$, the weighted centroid $\bmu(\cdot;\Uc)$ is $\Cc^\ell$ at~$\bx^*$.
\end{theorem}
\begin{proof}
    Suppose $i^+\in\Ic^+$ and $i^-\in\Ic^-$ are indices corresponding to mirrored constraints. Then using Lemma~\ref{lem:mirrored_halfspaces}, the region of integration is equivalently given by:
    \begin{align*}
    \Uc&_{\Ac\cup\Ic^-}(\bx)=\bigcap_{i\in\Ac\cup\Ic^-}\Omega_i(\bx)\\
    &= \Big(\bigcap_{i\in(\Ac\cup\Ic^-) \setminus\{i^+,i^-\}}\Omega_i(\bx)\Big) \cap \big(\Omega_{i^+}(\bx)\cap \Omega_{i^-}(\bx)\big)\\
    &= \Big(\bigcap_{i\in(\Ac\cup\Ic^-) \setminus\{i^+,i^-\}}\Omega_i(\bx)\Big) \cap \big(\Omega_{i^+}(\bx)\setminus \Omega_{i^-}^c(\bx)\big)\\
    &= \Big(\!\!\!\!\!\bigcap_{i\in(\Ac\cup\Ic^-) \setminus\{i^-\}}\!\!\!\!\!\Omega_i(\bx)\Big) \setminus  \Big(\!\!\!\!\!\bigcap_{i\in(\Ac\cup\Ic^-)\setminus\{i^+,i^-\}}\!\!\!\!\!\Omega_i(\bx)\cap \Omega_{i^-}^c(\bx)\Big).
    \end{align*}
    This is now expressed as a set difference between two regions, each of which involves one fewer mirrored constraint.  Thus, the integration over
    $\Uc_{\Ac\cup\Ic^-}(\bx)$ can be expressed as the difference between integrals over these two regions. The first region in the decomposition retains the index $i^+$ while eliminating its mirrored counterpart $i^-$, whereas the second region eliminates $i^+$ and instead introduces the complement halfspace $\Omega_{i^-}^c(\bx)$. Note that the constraint normal associated with the complement halfspace $\Omega_{i^-}^c(\bx)$ is $-\bb_{i^-}(\bx)$,
    which does not affect linear independence with respect to the remaining constraints.

    Repeating the decoupling process for each mirrored pair yields a finite collection of regions of integration, each of which is defined solely by constraints with linearly independent normals. The corresponding integrals are smooth with respect to $\bx$, as guaranteed by Proposition~\ref{prop:lin_ind}.    Therefore, it follows that the corresponding weighted centroid is smooth.
\end{proof}

Theorem~\ref{thm:ind_mirror} strengthens the earlier smoothness result obtained under strict linear independence in Proposition~\ref{prop:lin_ind}. It shows that smoothness of the weighted centroid can still be guaranteed even when some constraints are linearly dependent, provided that each dependent constraint appears as the unique mirrored counterpart of another constraint.

\begin{example}
[Triweight Formula for Input Bounds]\label{ex:input_bound}
    Mirrored structures arise naturally when safety requires independent outputs to be bounded from both sides. However, when input bounds, such as box constraints, are introduced alongside these state-dependent constraints, linear independence among all constraints considered no longer holds.
    We now
    illustrate the flexibility of our framework to handle state-independent input bounds through the choice of the weight function~$p$, enhancing the applicability of the results discussed in this section.

    For simplicity, consider a single-input ($m=1$) system with a mirrored pair of constraints:
    $$
    a_1(\bx) + b_1(\bx)u \leq 0, \qquad a_2(\bx) - b_1(\bx)u \leq 0,
    $$
    alongside an input bound $|u| \leq 1$. The collection of all four constraints does not meet the mirrored structure laid out in this section. Nevertheless, we may encode the state-independent constraint $|u| \leq 1$ within the choice of the weight function, cf. Remark~\ref{rmk:support}. In particular, we use here the \textit{triweight kernel}:
    $$
    p_\tri(u) = \begin{cases}
        \frac{35}{32}(1-u^2)^3, & |u|\leq1\\
        0, &|u|>1
    \end{cases}
    $$
    which is twice differentiable. Since the weighted centroid is computed over the support of $p_\tri$, the input bound is automatically enforced without any additional constraint. As a result, the admissible control map reduces to the mirrored pair alone, and smoothness of the weighted centroid follows from Theorem~\ref{thm:ind_mirror}. We obtain the following  formula:
    \begin{align*}
    \mu_\tri(\bx) &= \frac{\frac{u^2}{2}-\frac{3u^4}{4}+\frac{u^6}{2}-\frac{u^8}{8}\Big|_{u=\underline{c}(\bx)}^{u=\overline{d}(\bx)}}{u-u^3+\frac{3u^5}{5}-\frac{u^7}{7}\Big|_{u=\underline{c}(\bx)}^{u=\overline{d}(\bx)}},\\
    \underline{c}(\bx) &= \max\left(\min\Big(
    \frac{-a_1(\bx)}{b_1(\bx)},\frac{a_2(\bx)}{b_1(\bx)}\Big), -1\right),\\
    \bar{d}(\bx) &= \min\left(\max\Big(
    \frac{-a_1(\bx)}{b_1(\bx)},\frac{a_2(\bx)}{b_1(\bx)}\Big), 1\right).
    \end{align*}
    Despite the
    $\min$ and $\max$, the formula is twice differentiable if the collection of all constraints is strictly compatible.\hfill~$\bullet$
\end{example}

\subsection{Continuity when constraint normals vanish}
We now consider the case in which linear independence fails at a state $\bx^*$ because one or more constraint normals $\bb_i(\bx^*)$'s vanish. In this situation, the local coordinate transformation strategy used in the previous sections to render the domain of integration state-independent
is no longer applicable, as the set of linearly independent $\bb_i(\bx^*)$'s is not uniform in any neighborhood of such points. As a result, the smoothness analysis of the weighted centroid cannot proceed directly from the smoothness of the integrand. Nevertheless, the weighted centroid is continuous at such points under mild assumptions.

\begin{proposition}\longthmtitle{Continuity of Integrals when normals are zero}
\label{prop:continuous}
    Consider system~\eqref{sys:ctrl-affine} subject to admissible control maps $\{\Omega_i\}_{i\in\until{q}}$ in~\eqref{eq:ctrl-affine-feas}. Given
    $\bx^*$, define $\Ac^*=\setdefb{i\in\until{q}}{\bb_i(\bx^*)\neq \bzero}$ as the set of indices with nonzero $\bb_i(\bx^*)$. Assume that:
    \begin{enumerate}
        \item $a_i$, $\bb_i$, $\forall i\in \until{q}$ and $p$ are continuous at $\bx^*$;
        \item there exist integrable functions $\map{\Xi_0,\Xi_1}{\real^m}{\realnonneg}$, depending only on $\bu$, such that $|p(\bx,\bu)|\leq \Xi_0(\bu)$ and $\|\bu p(\bx,\bu)\|\leq \Xi_1(\bu)$ uniformly on a neighborhood of $\bx^*$;
        \item the normals $\{\bb_i(\bx^*)\}_{i\in{\Ac^*}}$ are linearly independent;
        \item the maps $\{\Omega_i\}_{i=1}^q$ are strictly compatible at $\bx^*$.
    \end{enumerate}
    Then  $\bmu(\bx^*;\Uc)=\bmu(\bx^*;\Uc_{\Ac^*})$, and the weighted centroid $\bmu(\cdot;\Uc)$ is continuous at $\bx^*$.
\end{proposition}
\begin{proof}
    We first show that $\Uc(\bx^*)=\Uc_{\Ac^*}(\bx^*)$.
    From the definition of $\Ac^*$, we have $\bb_i(\bx^*)=0$ for any index $i\in\until{q}\setminus\Ac^*$. By strict compatibility at $\bx^*$, it follows that for each $i\in\until{q}\setminus\Ac^*$:
    $$
    a_i(\bx^*)+\bb_i(\bx^*)\bu=a_i(\bx^*)<0,\quad\forall \bu\in\real^m.
    $$
    Therefore, the constraint indexed by $i\in\until{q}\setminus\Ac^*$ is satisfied trivially and $\Omega_i(\bx^*) = \real^m$. Consequently:
    \begin{align*}
        \Uc(\bx^*) = \bigcap_{i\in\until{q}}\Omega_i(\bx^*) = \bigcap_{i\in\Ac^*}\Omega_i(\bx^*) = \Uc_{\Ac^*}(\bx^*),
    \end{align*}
    which proves $\bmu(\bx^*;\Uc)=\bmu(\bx^*;\Uc_{\Ac^*})$.

    To establish continuity of the weighted centroid, it remains to show that the limit $\lim_{\bx\rightarrow \bx^*}\bmu(\bx;\Uc)$ exists and coincides with the value $\bmu(\bx^*;\Uc_{\Ac^*})$. For compactness of presentation, we first focus on the denominator of the weighted centroid~\eqref{eq:wc} and, for any arbitrary set $\Ac\subseteq \until{q}$, define:
    $$
    \eta_\Ac(\bx) =\int_{\Uc_\Ac(\bx)} p(\bx,\bu)d\bu.
    $$
    Note that the denominator of \eqref{eq:wc} is precisely $\eta_{\until{q}}$.

    Towards establishing continuity of $\eta_{\until{q}}$, we first show that $\eta_{\Ac^*}$ is continuous at $\bx^*$. Because the constraint normals $\{\bb_i(\bx^*)\}_{i\in\Ac^*}$ are continuous and linearly independent at~$\bx^*$, there exists a small $\delta_1>0$ such that all $\{\bb_i(\bx)\}_{i\in\Ac^*}$ remain nonzero and linearly independent for all $\bx\in\Ball{\delta_1}{\bx^*}$. As a result, we may apply the corresponding smooth local coordinate transformation $\bT(\bx,\bx^*)$  as in \eqref{eq:w-coordinates} to map to the $\bw$-coordinates, under which the domain of integration $\Omega_0$ is the same for all $\bx\in\Ball{\delta_1}{\bx^*}$. Since the integrand $p$ is continuous, it follows that $\eta_{\Ac^*}$ is continuous at~$\bx^*$.

    Next, we bound the difference:
    \begin{align*}
        |\eta_{\Ac^*}(\bx)-\eta_{\until{q}}(\bx)|
        &\leq \int_{\Uc_{\Ac^*}(\bx) \setminus \Uc(\bx)} |p(\bx,\bu)| d\bu\\
        &\leq \int_{\Uc_{\Ac^*}(\bx) \setminus \Uc(\bx)} \Xi_0(\bu) d\bu
    \end{align*}
    where the last inequality holds for all $\bx\in\Ball{\delta_\Xi}{\bx^*}$ with small $\delta_\Xi>0$. We will show  this term approaches zero as $\bx\rightarrow \bx^*$.

    For any $\bx$, the set difference in the domain of integration can be conservatively bounded as:
    \begin{align*}
    & \Uc_{\Ac^*}(\bx) \setminus \Uc(\bx) =
    \Uc_{\Ac^*}(\bx)\cap \Uc(\bx)^c   =\Uc_{\Ac^*}(\bx) \cap \Big(\bigcup_{j\in\until{q}}\!\Omega_j^c(\bx)\Big)\\
    & \quad =\Uc_{\Ac^*}(\bx) \cap \Big(\bigcup_{j\in\until{q}\setminus\Ac^*}\Omega_j^c(\bx)\Big) \subseteq\bigcup_{j\in\until{q}\setminus\Ac^*}\Omega_j^c(\bx).
    \end{align*}
    Here, for each $j$, the complement set can be bounded as:
    \begin{align*}
        \Omega_j^c(\bx)&=\setdefB{\bu\in\real^m}{a_j(\bx)+\bb_j(\bx)\bu > 0}\\
    &\subseteq\setdefB{\bu\in\real^m}{\|\bb_j(\bx)\|\|\bu\| > -a_j(\bx)}.
    \end{align*}
    Since $\bb_j(\bx^*)=0$ and $\bb_j$ is continuous for each $j\in\until{q} \setminus \Ac^*$, we have that for any $\epsilon_b>0$, there exists $\delta_2>0$ such that $\|\bb_j(\bx)\|<\epsilon_b$ for all $\bx\in\Ball{\delta_2}{\bx^*}$ and all $j\in\until{q} \setminus \Ac^*$.

    We then introduce a common lower bound on $-a_j(\bx)$ for all $j\in\until{q}\setminus \Ac^*$:
    $$
    \bar a =\min_{j\in\until{q}\setminus\Ac^*}\inf_{\bx\in\Ball{\delta_2}{\bx^*}}-a_j(\bx).
    $$
    Note that because strict compatibility at $\bx^*$ implies $a_j(\bx^*)<0$ for all $j\in\until{q}\setminus\Ac^*$, we can pick $\delta_2$ sufficiently small so that $\bar a$ is positive. Consequently, we derive:
    \begin{align*}
        \Omega_j^c(\bx)\subseteq\setdefB{\bu\in\real^m}{\|\bu\| > \bar a / \epsilon_b} ,
    \end{align*}
    for all $\bx\in\Ball{\delta_2}{\bx^*}$ and all $j\in\until{q} \setminus \Ac^*$. Since the superset is common for all $j\in\until{q} \setminus \Ac^*$, their union is also contained within the superset:
    \begin{align*}
        \bigcup_{j\in\until{q} \setminus \Ac^*}\Omega_j^c(\bx)\subseteq\setdefB{\bu\in\real^m}{\|\bu\| > \bar a / \epsilon_b}.
    \end{align*}
    Thus, we further bound the difference:
    $$
    |\eta_{\Ac^*}(\bx)-\eta_{\until{q}}(\bx)| \leq \int_{\{\bu\in\real^m\mid\|\bu\|\geq \bar a/\epsilon_b\}} \Xi_0(\bu)d\bu.
    $$
    Since $\Xi_0$ is integrable and non-negative, the integral evaluates to smaller values as $\epsilon_b$ decreases, independent of the state $\bx$.

    Given any $\epsilon>0$, there exists $\epsilon_b$ and a corresponding $\delta:= \min\{\delta_2,\delta_\Xi\}$ such that we can bound the difference $|\eta_{\Ac^*}(\bx)-\eta_{\until{q}}(\bx)|< \epsilon$ for all $\bx\in\Ball{\delta}{\bx^*}$. As a result, the function $\eta_{\Ac^*}-\eta_{\until{q}}$ converges to zero at $\bx^*$. Since $\eta_{\Ac^*}$ is continuous at $\bx^*$, we derive $\lim_{\bx\rightarrow\bx^*}\eta_{\until{q}}(\bx)=\eta_{\Ac^*}(\bx^*)$. An analogous argument applies to the numerator of~\eqref{eq:wc}, using $\Xi_1$ for each component of the integrand $u_i p(\bx,\bu)$. Therefore, the limit $\lim_{\bx\rightarrow\bx^*}\bmu(\bx;\Uc)$ exists and coincides with the value $\bmu(\bx;\Uc_{\Ac^*})$, establishing that the weighted centroid $\bmu(\cdot;\Uc)$ is continuous at $\bx^*$, concluding the proof.
\end{proof}

Proposition~\ref{prop:continuous} establishes continuity of the weighted centroid at states where some $\bb_i$'s vanish, under an additional assumption on the uniform integrability of the weight function. This assumption ensures the tails of the integrands are uniformly controlled in a neighborhood of $\bx^*$. For a general
choice of $p$, this assumption may not hold. Nevertheless, since the weight function is a design choice, this assumption is not restrictive in practice. For instance, we can select weights that are simply independent of the state, or adopt structured families such as Gaussian weights with continuously varying covariance like~\eqref{eq:gaussian_one}, which satisfy the assumption. In the next section, we show that Gaussian weights admit additional structure that enables stronger smoothness results.

\subsection{Gaussian Weights: Beyond Continuity}
This section investigates the use of Gaussian weights to handle the case with two ($q=2$) admissible control maps. Throughout this section, we assume the feasible set $\Omega_1(\bx)\cap\Omega_2(\bx)$ is nonempty for all $\bx$ under consideration. Consider the multivariate Gaussian distribution with a diagonal covariance matrix $\Sigma(\bx)=\sigma(\bx)^2\bI_{m}$, where $\sigma(\bx)>0$:
\begin{equation}\label{eq:gaussian_gen}
p_\gau^m = ((2\pi)^{-\frac{1}{2}} \sigma(\bx))^m\exp\left(-\frac{\bu^\top\bu}{2\sigma(\bx)^2}\right).
\end{equation}
When the number of admissible control maps is two ($q=2$), the weighted centroid~\eqref{eq:wc} using this $p_\gau^m$ as its weight function:
$$
\bmu_\gau(\bx;\Omega_1\cap\Omega_2) = \frac{\int_{\Omega_1(\bx)\cap \Omega_2(\bx)}\bu p_\gau^m(\bx,\bu)d\bu}{\int_{\Omega_1(\bx)\cap \Omega_2(\bx)}p_\gau^m(\bx,\bu)d\bu}
$$
admits a closed-form solution (see also~\cite{GMT:65}).

For
brevity,
define the standardized margins $\zeta_i(\bx)= -a_i(\bx)/(\sigma(\bx)\|\bb_i(\bx)\|)$
and the correlation $\rho(\bx)=\bb_1(\bx)^\top\bb_2(\bx)/(\|\bb_1(\bx)\|\|\bb_2(\bx)\|)$ between the two constraint normals $\bb_1$ and $\bb_2$. In addition, let $(\zeta_1,\zeta_2)\mapsto \Phi^2(\zeta_1,\zeta_2;\bG(\bx))$ denote the cumulative distribution function (CDF) of the zero-mean bivariate normal distribution with covariance matrix:
$$
\bG(\bx)=\begin{bmatrix}
    1 & \rho(\bx) \\
    \rho(\bx) & 1
\end{bmatrix}.
$$
With this setup, the Gaussian-weighted centroid is given by:
\begin{multline}\label{eq:closed_form}
\bmu_\gau(\bx;\Omega_1\cap\Omega_2) =\\
\begin{cases}
\bzero, &~\bb_1(\bx) = 0, \bb_2(\bx) = 0\\
\bmu_\gau(\bx;\Omega_1), &~\bb_1(\bx)\neq 0, \bb_2(\bx) = 0 \\
\bmu_\gau(\bx;\Omega_2), &~\bb_1(\bx) = 0, \bb_2(\bx) \neq 0 \\
\bmu_{12}(\bx;\Omega_1\cap\Omega_2), &~\bb_1(\bx) \neq 0, \bb_2(\bx)\neq 0
\end{cases}.
\end{multline}
where:
\begin{align*}
&\bmu_\gau(\bx;\Omega_i) \defeq -\sigma \frac{\bb_i^\top}{\|\bb_i\|} \frac{p_\gau^1(\zeta_i)}{\Phi(\zeta_i)},\\
&\bmu_{12}(\bx;\Omega_1\cap\Omega_2) \defeq\\
&\begin{cases}-\frac{\sigma}{\Phi^2(\zeta_1,\zeta_2; \bG)}\Big(
    \frac{\bb_1^\top}{\|\bb_1\|}p_\gau^1(\zeta_1)\Phi\big(\frac{\zeta_2-\rho\zeta_1}{\sqrt{1-\rho^2}}\big) \\\qquad+ \frac{\bb_2^\top}{\|\bb_2\|}p_\gau^1(\zeta_2)\Phi\big(\frac{\zeta_1-\rho\zeta_2}{\sqrt{1-\rho^2}}\big) \Big),
& \rho \neq \pm 1 \\
\sigma \frac{\bb_1^\top}{\|\bb_1\|} \frac{p_\gau^1(-\zeta_2)-p_\gau^1(\zeta_1)}{\Phi(\zeta_1)-\Phi(-\zeta_2)},    &\rho = -1 \\
 -\sigma \frac{\bb_1^\top}{\|\bb_1\|}\max\left\{\frac{p^1_\gau\left(\zeta_1\right)}{\Phi\left(\zeta_1\right)},\frac{p^1_\gau\left(\zeta_2\right)}{\Phi\left(\zeta_2\right)}\right\},   &\rho = 1
\end{cases}
\end{align*}
Here, for brevity, we have omitted the dependency on $\bx$ for the functions $a_i$, $\bb_i$, $\zeta_i$, $\sigma$, $\rho$, and $\bG$ used above.
The closed-form expressions~\eqref{eq:closed_form} above allow us to directly characterize the regularity of the weighted centroid under Gaussian weights.

\begin{proposition}\label{prop:gau}
\longthmtitle{Smoothness of Gaussian-weighted centroid}
    Consider the weighted centroid function~\eqref{eq:wc} associated with two admissible control maps $\{\Omega_1,\Omega_2\}$ of the form~\eqref{eq:ctrl-affine-feas}, under the Gaussian weight function~\eqref{eq:gaussian_gen}. Assume that:
    \begin{enumerate}
        \item $a_i$, $\bb_i$, and $\sigma$ are $\Cc^\ell$;
        \item the maps $\{\Omega_1,\Omega_2\}$ are strictly compatible at each $\bx$.
    \end{enumerate}
    Then the Gaussian-weighted centroid is $\Cc^\ell$ for all $\bx$ such that $\rho(\bx)\neq 1$.
\end{proposition}
\begin{proof}
The claim follows by inspection of
\eqref{eq:closed_form}. For each case where $\rho(\bx)\neq 1$, the expression is a composition of smooth functions of $a_i(\bx)$, $\bb_i(\bx)$, and $\sigma(\bx)$, together with the Gaussian PDF and CDF, and is therefore smooth in $\bx$ on the corresponding region. In addition, the limiting behavior as $\rho(\bx)\to -1$ and as $\bb_i(\bx)\to\bzero$ matches the corresponding case in~\eqref{eq:closed_form}, yielding smooth transitions across these boundaries. The exponential decay of the Gaussian ensures that all derivatives of the closed-form expression remain bounded in the limits $\bb_i(\bx)\to\bzero$ and $\rho(\bx)\to -1$, since derivatives introduce at most polynomial growth in $\zeta_i(\bx)$, which is dominated by the Gaussian exponential term (cf. Theorem~\ref{thm:Leibniz_converging}).
\end{proof}

The only source of nonsmoothness in the Gaussian weighted centroid $\bmu_\gau$ arises in the case $\rho(\bx)=1$, due to the presence of the maximum operator. To eliminate this nonsmoothness, we introduce a smooth surrogate defined as the sum of the single-constraint Gaussian centroids:
\begin{equation}\label{eq:smooth-rho=1}
\bmu_{\rm sum}(\bx;\Omega_1\cap \Omega_2) \defeq \bmu_\gau(\bx;\Omega_1)+\bmu_\gau(\bx;\Omega_2)
\end{equation}
By construction, this function is smooth for all $\bx$, so long as the constraints satisfy strict compatibility. However, it does not necessarily satisfy both constraints. The following result shows that it is admissible whenever $\rho(\bx)\geq 0$.

\begin{lemma}\label{lem:smooth_select}
\longthmtitle{Admissibility of $\bmu_{\rm sum}$}
    Consider two admissible control maps $\{\Omega_1,\Omega_2\}$ of the form~\eqref{eq:ctrl-affine-feas}, and let $\bmu_{\rm sum}$ be defined by~\eqref{eq:smooth-rho=1}. Then $\bmu_{\rm sum}(\bx)\in \Omega_1(\bx)\cap\Omega_2(\bx)$ for each $\bx$ such that $\rho(\bx)\geq 0$.
\end{lemma}
\begin{proof}
Since $\rho(\bx)\geq 0$, the products $\bb_1(\bx)\bmu_\gau(\bx;\Omega_2) \leq 0$ and $\bb_2(\bx)\bmu_\gau(\bx;\Omega_1) \leq 0$ are nonpositive. Thus, we deduce for each $i\in\{1,2\}$,
    \begin{align*}
    a_i(\bx) +\bb_i(\bx)\bmu_{\rm sum}(\bx;\Omega_1\cap \Omega_2) &\leq  a_i(\bx) +\bb_i(\bx)\bmu_\gau(\bx;\Omega_i).
    \end{align*}
Since $\bmu_\gau(\bx;\Omega_i)\in\Omega_i(\bx)$ (because it is a weighted centroid of the set~\cite[p. 25]{SB-LV:09}), we have:
$$
a_i(\bx) +\bb_i(\bx)\bmu_\gau(\bx;\Omega_i) \leq 0.
$$
Hence, we conclude $\bmu_{\rm sum}(\bx)\in \Omega_i(\bx)$ for~$i\in\{1,2\}$.
\end{proof}

Since the surrogate $\bmu_{\rm sum}$ is admissible whenever $\rho(\bx)\ge 0$, we  construct a globally smooth Gaussian-weighted centroid by smoothly transitioning from $\bmu_\gau$ to $\bmu_{\rm sum}$  over $\rho(\bx)\in[0,1]$. Let $s:\real\to[0,1]$ be a $C^\infty$ transition function satisfying:
$$
s(r)=0 \ \text{for } r\le 0,
\qquad
s(r)=1 \ \text{for } r\ge 1,
$$
One example of this is given by:
\begin{equation}\label{eq:smooth_transition}
  s(r) = \begin{cases}0, & r\leq 0 \\
    \left(1+\frac{e^{1/r}}{e^{1/(1-r)}}\right)^{-1}, & 0<r<1 \\
    1,& r\geq 1
  \end{cases}.
\end{equation}
We then define the selection function:
\begin{equation}\label{eq:selection_function}
\bk(\bx) \defeq  \big(1-s(\rho(\bx))\big)\bmu_\gau(\bx;\Omega_1\cap\Omega_2)+ s(\rho(\bx))\bmu_{\rm sum}(\bx).
\end{equation}
The following result shows that $\bk$ is smooth for all $\bx$, and satisfies both constraints.

\begin{proposition}\longthmtitle{Selection function via transition}\label{prop:two_constraints}
    Consider two admissible control maps $\{\Omega_1,\Omega_2\}$ of the form~\eqref{eq:ctrl-affine-feas}. Let the selection function~$\bk:\real^n\to\real^m$ be defined by~\eqref{eq:selection_function}. The function $\bk$ satisfies both constraints at each $\bx$, i.e., $\bk(\bx)\in\Omega_1(\bx)\cap\Omega_2(\bx)$. Moreover, assume that:
    \begin{enumerate}
        \item $a_i$, $\bb_i$, and $\sigma$ are $\Cc^\ell$;
        \item the maps $\{\Omega_1,\Omega_2\}$ are strictly compatible at each $\bx$.
    \end{enumerate}
    Then the selection function~$\bk$ is $\Cc^\ell$.
\end{proposition}
\begin{proof}
    If $\rho(\bx)<0$, then $\bk(\bx)=\bmu_\gau(\bx;\Omega_1\cap\Omega_2)$, which lies in $\Omega_1(\bx)\cap\Omega_2(\bx)$ since it is a weighted centroid of that convex set~\cite[p. 25]{SB-LV:09}. If $\rho(\bx)\geq 0$, then $\bk(\bx)$ is a convex combination of $\bmu_\gau(\bx;\Omega_1\cap\Omega_2)$ and $\bmu_{\rm sum}(\bx)$, both of which belong to the intersection $\Omega_1(\bx)\cap\Omega_2(\bx)$ (cf. Lemma~\ref{lem:smooth_select} for $\bmu_{\rm sum}$). Since the intersection $\Omega_1(\bx)\cap\Omega_2(\bx)$ is convex, it follows that $\bk(\bx)$ must lie in the set.    Smoothness follows from Proposition~\ref{prop:gau}, the smoothness of the transition function $s$, and the smoothness of $\bmu_{\rm sum}$ under strict compatibility ($\zeta_i\rightarrow+\infty$ as $\bb_i\rightarrow \bzero$).
\end{proof}
By introducing the surrogate $\bmu_{\rm sum}$ and blending it with the original Gaussian-weighted centroid, we obtain a globally smooth selection function $\bk$ that remains admissible for all $\bx$. This construction preserves the geometric interpretation of the weighted centroid while eliminating nonsmooth behavior.

\section{Corollaries for Stabilization and Safety}\label{sec:control}
The previous sections presented general results for constructing smooth feedback controllers enforcing classes of affine input constraints. We now specialize these results to the setting of stabilization and safety with CLFs and CBFs. In particular, this section provides the precise statements on how the proposed weighted centroid method yields smooth stabilizing and/or safeguarding controller designs.

\subsection{Safety Filtering with CBFs}
We begin by considering safety filtering using CBFs. The objective is typically to modify a given nominal controller $\bk_{\rm{d}}(\bx)$ while ensuring that the resulting control input satisfies the CBF constraint \eqref{eq:cbf_deriv}, which is a sufficient condition for forward invariance of the safe set $\mathcal{S}$ defined in \eqref{eq:safeset}.

The following result applies the weighted centroid methodology for a single affine constraint, as developed in Sec.~\ref{sec:single_ineq}.
\begin{corollary}\longthmtitle{Safeguarding Formulas}
\label{cor:safe_formulas}
    Let $h$ be a CBF for~\eqref{sys:ctrl-affine} on a set $\Sc$ as in~\eqref{eq:safeset}. Given a $\Cc^\ell$ nominal controller $\map{\bk_{\rm{d}}}{\real^n}{\real^m}$, let
    $a_1(\bx)=-\Lie_\bf h(\bx) - \Lie_{\bg}h(\bx)\bk_{\rm{d}}(\bx) - \alpha_h(h(\bx))$ and $\bb_1(\bx)=-\Lie_\bg h(\bx)$. If the assumptions on $a_1$, $\bb_1$, and $p$, in Theorem~\ref{thm:Leibniz_converging} hold, then the controller $\map{\bk}{\real^n}{\real^m}$ defined  in \eqref{eq:single-input-controller-and-map} yields a feedback law $\bu= \bk_{\rm{d}}(\bx) + \bk(\bx)$ that is $\Cc^\ell$ on $\Sc$ and renders $\Sc$ forward invariant.
\end{corollary}
\begin{proof}
Define the coordinate transformation with $\bu=\bk_{\rm{d}}(\bx)+\tilde\bu$ where $\tilde \bu\in\real^m$ represents the control adjustment to $\bk_{\rm{d}}(\bx)$. Consider the nonlinear system given by:
$$
\dot \bx = \bf(\bx)+\bg(\bx)\bk_{\rm{d}}(\bx)+\bg(\bx)\tilde\bu.
$$
The forward invariance result follows from applying Theorem~\ref{thm:Leibniz_converging} to this equivalent system.
\end{proof}

This result recovers a smooth safety filter for a single CBF constraint. Alternative solutions for this problem are well documented
\cite{MHC-PO-GB-ADA:23-csl} and can, for instance, be obtained via Sontag-type formulas, see Example~\ref{ex:sontag}.
The strength of the weighted centroid methodology is more apparent when multiple constraints are present. To this end, we now consider multiple mirrored constraints with linearly independent normals.

\begin{corollary}\longthmtitle{Linearly Independent Mirrored CBFs}\label{cor:li-m-CBFs}
Let $\{h_i\}_{i=1}^q$ be CBFs for~\eqref{sys:ctrl-affine} defining safe sets $\{\Sc_i\}_{i=1}^q$ as in~\eqref{eq:safeset}. Given a $\Cc^\ell$ nominal controller $\map{\bk_{\rm{d}}}{\real^n}{\real^m}$, define the function $a_i(\bx)=-\Lie_\bf h_i(\bx) - \Lie_{\bg}h_i(\bx)\bk_{\rm{d}}(\bx) - \alpha_{h_i}(h_i(\bx))$ and $\bb_i(\bx)=-\Lie_\bg h_i(\bx)$, for each $i\in\until{q}$. If the assumptions on $a_i$, $\bb_i$, and $p$, in Theorem~\ref{thm:ind_mirror} hold, then the controller~$\bmu(\cdot;\cap_{i=1}^q \Omega_i(\bx))$ based on the weighted centroid of the admissible control map $\{\Omega_i\}_{i=1}^q$ defined  in \eqref{eq:ctrl-affine-feas} yields a feedback law $\bu= \bk_{\rm{d}}(\bx) + \bmu(\bx;\cap_{i=1}^q \Omega_i(\bx))$ that is $\Cc^\ell$ on $\Sc=\cap_{i\in\until q}\Sc_i$ and renders $\Sc$ forward invariant.~\hfill$\blacksquare$
\end{corollary}

Corollary~\ref{cor:li-m-CBFs} applies to collections of mirrored CBF constraints with linearly independent normals, while accommodating the dependent case of mirrored constraints. Although this assumption may appear restrictive, it is often satisfied in practice. In particular, the commonly used optimal-decay formulation~\cite{JZ-BZ-ZL-KS:21,PO-MHC-TGM-ADA:25-cdc} introduces auxiliary variables that render the resulting constraint normals linearly independent for states in the interior of the safe set. Specifically, instead of enforcing the standard CBF constraint in~\eqref{eq:cbf_deriv}, the formulation uses:
$$
\Lie_\bf h(\bx)+\Lie_\bg h(\bx)\bu > -\omega \alpha_h(h(\bx))
$$
where $\omega\in\realnonneg$ is an auxiliary control variable. This reformulation yields linear independence of the CBF constraint normals but introduces the additional constraint $\omega \geq 0$. Nevertheless, this additional constraint is state-independent and can be handled within our framework through an appropriate choice of the weighting function $p$'s support, cf. Example~\ref{ex:input_bound}.

\subsection{Stabilization with CLFs}
We now turn to the problem of stabilization using CLFs. The weighted
centroid methodology yields a smooth stabilizing controller on
$\R^n \setminus \{\bzero\}$, as stated in the following result.
\begin{corollary}\longthmtitle{Stabilizing Formulas}
\label{cor:stab_formulas}
   Let $V$ be a CLF for~\eqref{sys:ctrl-affine} and define the admissible control map $\Omega_v$ as in \eqref{eq:single-input-controller-and-map} with $a_1(\bx)=\Lie_\bf V(\bx) + \alpha(V(\bx))$ and $\bb_1(\bx)=\Lie_\bg V(\bx) $, and let $p$ be any weight function satisfying the hypotheses of Theorem \ref{thm:Leibniz_converging}. Provided $a_1,\bb_1\in\Cc^l$, then the controller $\bu=\bk(\bx)$, with $\bk$ as in \eqref{eq:single-input-controller-and-map}, is $\Cc^l$ on $\R^n\setminus\{\bzero\}$ and renders the origin globally asymptotically stable.~\hfill$\blacksquare$
\end{corollary}

Corollary~\ref{cor:stab_formulas} establishes asymptotic stability using the CLF constraint on $\real^n\setminus\{0\}$ together with the convention $\bk(\bzero)=\bzero$. In contrast to the CBF setting, applying the weighted centroid methodology to CLF-based stabilization requires additional care near the origin $\bx=\bzero$. The
issue arises from the fact that \emph{strict} feasibility of the CLF constraint~\eqref{eq:clf_deriv} fails at the origin. Indeed, since $\Lie_{\bf}V(\bzero)=0$ and $\Lie_{\bg}V(\bzero)=\bzero$, the inequality cannot be satisfied strictly at $\bx=\bzero$. Thus, the weighted centroid $\bmu$, while smooth on $\real^n \setminus \{\bzero\}$, may fail to be continuous at the origin.
This mirrors a classical result: Sontag's universal formula, the canonical smooth stabilizing controller, is also only guaranteed to be continuous at the origin under the \emph{small control property}.

\begin{definition}\longthmtitle{Small Control Property}
{\rm
  A CLF $V$ satisfies the \emph{small control property} if for every $\epsilon>0$, there exists
  $\delta>0$ such that if $\bx\neq\bzero$ and $\|\bx\|<\delta$, there exists $\bu\in
  \R^m$ with $\|\bu\|<\epsilon$ such that:
  $$\Lie_{\bf}V(\bx) + \Lie_{\bg}V(\bx)\bu < - \alpha(V(\bx)).$$
  }
\end{definition}

The small control property is in fact necessary for the existence of any continuous controller satisfying the CLF constraint locally at the origin, as any such controller $\bk$ with $\bk(\bzero)= \bzero$ must take arbitrarily small values near the origin while still satisfying the CLF constraint.

To leverage the small control property, our approach involves modifying the weight function by introducing a state-dependent scaling on the control input $\bu$. The idea is motivated by the Gaussian weight $p_{\rm gau}^1$ in~\eqref{eq:gaussian_one} with a smooth positive definite state-dependent standard deviation $\map{\sigma}{\real^n}{\realnonneg}$ satisfying $\sigma(\bzero) = 0$ and $\sigma(\bx)>0$ elsewhere. As $\bx \to \bzero$, the standard deviation $\sigma(\bx) \to 0$ concentrates the Gaussian weight near $\bu = \bzero$, pulling the weighted centroid toward the origin. Effectively, the weighted centroid approaches a min-norm selection near the origin, which is known to be continuous under the small control property, as shown below. This intuition motivates the conditions in the following result that enables continuity at the origin.

\begin{theorem}\longthmtitle{Continuity at the Origin}\label{thm:scp}
Let $V$ be a CLF for a single-input ($m=1$) control-affine system~\eqref{sys:ctrl-affine}  satisfying the small control property. Define $a_1(\bx) =
\Lie_{\bf}V(\bx) + \alpha(V(\bx))$ and $b_1(\bx) =
\Lie_{\bg}V(\bx)$. Let $\map{p_0}{\real}{\realnonneg}$ be a $\Cc^1$ weight function such that $v \mapsto p_0(v)$ is integrable and $v \mapsto v\,p_0(v)$ is integrable with $\int_{-\infty}^{\infty} v\,p_0(v)\,dv = 0$. Let $\map{\sigma}{\real^n}{\realnonneg}$ be a
smooth positive definite function satisfying $\sigma(\bzero) = 0$.
Define the scaled weight $p(\bx, v) = p_0(v/\sigma(\bx))$
and let $\mu$ be the weighted centroid controller~\eqref{eq:wc}
under $p$ applied to the admissible control map $\Omega_1(\bx) =
\setdefb{v\in\real}{a_1(\bx) + b_1(\bx) v \leq 0}$. If $\Omega_1$ is strictly compatible on $\real^n\setminus \{\bzero\}$ and the following tail limits on the weight:
    \begin{equation}\label{eq:tail_condition}
\lim_{v\rightarrow \pm\infty} \tfrac{p_0(v)+v\frac{dp_0}{dv}(v)}{2p_0(v)+v \frac{dp_0}{dv}(v)} = R^\pm
\end{equation}
are finite for some $R^-,R^+>0$, then $\mu$ is continuous at $\bx = \bzero$ with $\mu(\bzero) = 0$.
\end{theorem}
\begin{proof}
    From strict compatibility, $\mu(\bx)=0$ when $b_1(\bx)=0$ and $\bx\neq \bzero$. For the remainder of this proof, we assume without loss of generality (relying on the corresponding limit $R^-$ or $R^+$) that $b_1(\bx) > 0$. Applying the substitution $w = v/\sigma(\bx)$, the weighted centroid simplifies to:
    $$
    \mu(\bx) = \sigma(\bx)\frac{\int_{-\infty}^{d^\sigma(\bx)} w  p_0(w) dw}{\int_{-\infty}^{d^\sigma(\bx)}  p_0(w) dw}\defeq \sigma(\bx) \mu^\sigma(d^\sigma(\bx))
    $$
    where $d^\sigma(\bx) = -a_1(\bx)/(\sigma(\bx)b_1(\bx))$ and $\map{\mu^\sigma}{\real}{\real}$ is the truncated mean of $p_0$ as a function of the upper limit $d^\sigma$.

    We further observe that $\mu^\sigma$ is monotonically increasing in $d^\sigma$, since:
    \begin{align*}
    \frac{d\mu^\sigma}{dd^\sigma} &= \frac{d^\sigma p_0(d^\sigma)\int_{-\infty}^{d^\sigma} p_0(w)\, dw - p_0(d^\sigma)\int_{-\infty}^{d^\sigma} w\,p_0(w)\, dw}   {\left(\int_{-\infty}^{d^\sigma} p_0(w)\, dw\right)^2} \\
    &= \frac{p_0(d^\sigma)}{\int_{-\infty}^{d^\sigma} p_0(w)\, dw}\left(d^\sigma - \mu^\sigma(d^\sigma)\right) \geq 0,
    \end{align*}
    where the inequality holds because $d^\sigma \geq \mu^\sigma(d^\sigma)$,
    as $\mu^\sigma(d^\sigma)$ is the weighted centroid of the interval $[-\infty, d^\sigma]$.

    From the small control property, if $\bx\in\Ball{\delta}{\bzero}$ for some small $\delta>0$, then there is $v^*\in\real$ with $|v^*| < \epsilon$ such that:
    $$
    a_1(\bx)+b_1(\bx)v^* \leq 0 \implies -\epsilon \leq v^* \leq -\frac{a_1(\bx)}{b_1(\bx)}.
    $$
    Therefore, $d^\sigma(\bx) \geq -\epsilon/\sigma(\bx)$, and from monotonicity:
    $$
      \mu^\sigma\big(\tfrac{-\epsilon}{\sigma(\bx)}\big) \leq \mu^\sigma(d^\sigma(\bx))\leq \mu^\sigma(\infty) =0.
    $$

    To study the behavior as $\bx\to\bzero$, we analyze the limit:
    \begin{align*}
    \lim_{d\to-\infty}\frac{\mu^\sigma(d)}{d}
    &= \lim_{d\to-\infty}\frac{\int_{-\infty}^{d} wp_0(w)\,dw}
    {d\int_{-\infty}^{d} p_0(w)\,dw}  .
    \end{align*}
    We note here that the numerator evaluates to zero in the limit, and the denominator also evaluates to zero because:
    \begin{align*}
    \lim_{d\to -\infty}\left| d\int_{-\infty}^{d} p_0(w)\,dw \right| &\leq \lim_{d\to -\infty} \int_{-\infty}^{d} |d| p_0(w)\,dw \\
    &\leq \lim_{d\to -\infty} \int_{-\infty}^{d} |w| p_0(w)\,dw = 0
    \end{align*}
    Therefore, we apply L'H\^{o}pital's rule to obtain:
    \begin{align*}
    \lim_{d\to-\infty}\frac{\mu^\sigma(d)}{d}&= \lim_{d\to-\infty}
    \frac{dp_0(d)}{\int_{-\infty}^{d} p_0(w)\,dw + dp_0(d)} \\
    &= \lim_{d\to-\infty}
     \frac{p_0(d)+d\frac{dp_0}{dv}(d)}{2p_0(d)+d \frac{dp_0}{dv}(d)} = R^-,
    \end{align*}
    where the second equality follows from a second application of L'H\^{o}pital's rule (Note $dp_0(d)\to 0$ from integrability of $vp_0$). Therefore, we deduce:
    $$
    -R^-\epsilon = \lim_{\bx\rightarrow \bzero}\sigma(\bx)\mu^\sigma\left(\frac{-\epsilon}{\sigma(\bx)}\right) \leq \mu(\bx) \leq 0,
    $$
    which holds for any $\epsilon > 0$. Taking $\epsilon \to 0$ gives $\lim_{\bx\to\bzero}\mu(\bx) = 0$, establishing continuity at the origin.
\end{proof}

Theorem~\ref{thm:scp} gives a condition on the weight $p$ such that the resulting weighted centroid is continuous under the small control property. The tail condition~\eqref{eq:tail_condition} may appear unusual, but it is straightforward to verify for a given weight function $p_0$ by direct computation. For instance, one may verify that the Gaussian weight $p_0(v) = e^{-v^2/2}$ yields $R^\pm = 1$, while Sontag's weight $p_0(v) = (v^2+1)^{-3/2}$ yields $R^\pm = 2$. In both cases, the condition is satisfied, confirming that the weighted centroid under these weights is continuous at the origin.

The result of Theorem~\ref{thm:scp} is stated for single-input systems, but it extends naturally to the multi-input case ($m>1$) with a single ($q=1$) constraint via the technique introduced in Theorem~\ref{thm:Leibniz_converging} to reduce the constraint into a single-input one.

\subsection{Safe Stabilizing with a CLF and a CBF}
Finally, we consider the simultaneous enforcement of stability and safety objectives using a CLF and a CBF. In this setting, it is unlikely linear independence will hold. The following result directly follows from Proposition~\ref{prop:two_constraints}.

\begin{corollary}\longthmtitle{Safe Stabilizing Formulas}
    \label{cor:safe_stab}
    Let $V$ be a CLF and $h$ be a CBF for~\eqref{sys:ctrl-affine} defining a safe set $\Sc$ as in~\eqref{eq:safeset}. Define admissible control maps $\{\Omega_1,\Omega_2\}$ as in Corollaries~\ref{cor:safe_formulas} and \ref{cor:stab_formulas}. Under the assumptions of Proposition~\ref{prop:two_constraints} and assuming that the origin is in the interior of $\Sc$, the controller $\bu=\bk(\bx)$, with $\bk$ defined in~\eqref{eq:selection_function} on $\Sc \setminus \{\bzero\}$, is $\Cc^\ell$ and renders the set $\Sc$ safe. Furthermore, if $\bk$ is extended to the origin by setting $\bk(\bzero)=\bzero$, then the origin is an equilibrium of the closed-loop system and is asymptotically stable.~\hfill$\blacksquare$
\end{corollary}

Corollary~\ref{cor:safe_stab} provides a safe-stabilizing formula for selecting control inputs satisfying both a CLF and a CBF constraint at each $\bx$ on $\Sc\setminus\{\bzero\}$. The condition $\bk(\bzero)=\bzero$ is stated explicitly because unlike in the CLF-only case, the controller $\bk$ in \eqref{eq:selection_function} may not be zero at the origin due to the presence of the CBF constraint. This condition ensures that the origin remains an equilibrium of the closed-loop system so that we may discuss stability. Note importantly that since $\bzero\in\interior(\Sc)$ is in the interior of the safe set, assigning $\bk(\bzero)=\bzero$ does not compromise safety.
As in the CLF-only case, continuity at the origin requires additional care because the strict compatibility assumption necessarily fails at $\bx=\bzero$. In the CLF-CBF setting, the first step is to ensure that $\bu=\bzero$ satisfies the CBF constraint. Since $\bzero\in\interior(\Sc)$, this can be achieved through picking a sufficiently large class-$\Kc$ function $\alpha_h$. Once $\bu=\bzero$ is feasible for the CBF constraint at $\bx=\bzero$, choosing $\sigma$ to be a positive definite function renders continuity at the origin, as evident from the closed-form formula in \eqref{eq:selection_function}.

The results in this section show how the weighted centroid methodology provides a systematic path from CLF- and CBF-induced admissible control maps to explicit smooth feedback formulas. By interpreting CLF and CBF constraints as affine constraints on the input, the framework uses a geometric selection rule that averages over the feasible set. This perspective not only recovers classical universal formulas in the single-constraint setting, but also provides a natural pathway toward smooth controller synthesis in multiple-constraint settings.

\section{Simulations and Experiments}\label{sec:experiments}
We present complementary simulation and hardware studies illustrating the efficacy of our weighted centroid methodology. Consider a quadrotor operated through a velocity-command interface, with a pair of safety constraints, $h_1 (\bx) =  x_1 + 1$ and $h_2(\bx) = 1 - x_1$.
 At the outer-loop level, we model the position dynamics as the single-integrator system $\dot x_1 = u$,
where $u$ denotes the velocity command sent to the quadrotor. The desired velocity command $k_\des$ is to track a desired position trajectory $r_\des$, via a tracking controller with a feedforward term:
$$
k_\des(x_1,t) = 2.5(r_\des(t)-x_1)+\dot r_\des(t).
$$

In both the simulation and experimental studies, we compare the standard CBF-QP safety filter in~\cite{ADA-SC-ME-GN-KS-PT:19} with the smooth safety filter obtained by applying the mirrored Sontag formula in~\eqref{eq:sontag_mirror} to the CBF constraints according to Corollary~\ref{cor:li-m-CBFs}\footnote{We use ${\alpha_h(r)=5r}$ in the simulation and ${\alpha_h(r)=3r}$ in the experiment.}. We denote the resulting filtered velocity command by $v_{\rm cmd}$.

\subsection{Simulation: Quadruple Integrator}
To investigate the effect of lower-level tracking dynamics in realizing $v_{\rm cmd}$, we explicitly model the higher-order closed-loop dynamics in simulation. Motivated by the fact that, after dynamic extension and feedback linearization, the dynamics of a quadrotor take the form of a quadruple integrator, we consider the following simplified closed-loop model:
\begin{align}\label{eq:quadruple-integrator}
    \dot x_1 &= x_2, \qquad \dot x_2 = x_3, \qquad \dot x_3 = x_4 \\
    \dot x_4 &= -k_1(x_2\!-\!v_{\rm cmd})-k_2(x_3\!-\!\dot v_{\rm cmd})-k_3(x_4\!-\!\ddot v_{\rm cmd}) \nonumber
\end{align}
Our simulation uses $k_1=3$, $k_2=6$, and $k_3=12$ to obtain relatively slow tracking dynamics, thereby accentuating the effect of variations in $v_{\rm cmd}$ on the lower-level tracking response. This choice is intended to illustrate why smoothness of $v_{\rm cmd}$ can be important when the filtered command must be realized through higher-order tracking dynamics.

\begin{figure}
    \centering
    \includegraphics{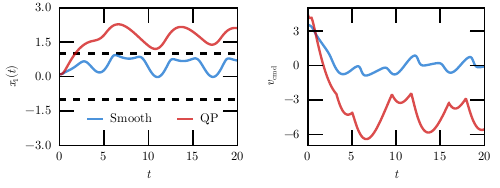}
    \caption{Simulation results of the quadruple integrator tracking a velocity command generated by our weighted centroid approach against a standard QP approach. \textbf{Left:} the position $x_1$ of the quadruple integrator when tracking the smooth (blue) vs QP (red) velocity, where the black dashed curves denote the state bounds. \textbf{Right:} corresponding filtered velocity commands. }
    \label{fig:states}
    \vspace{-10pt}
\end{figure}

\begin{figure}
    \centering
    \includegraphics{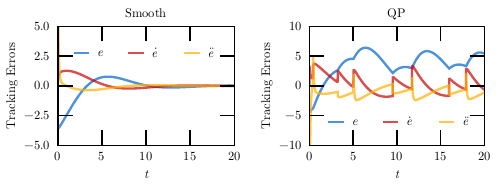}
    \caption{Comparison of the tracking errors in simulations when generating a safe reference velocity using the mirrored Sontag formula~\eqref{eq:sontag_mirror} (\textbf{left}) vs the CBF-QP approach (\textbf{right}).}
    \label{fig:errors}
    \vspace{-15pt}
\end{figure}

The simulation uses the following reference trajectory:
$$
r_d(t) = 1-\tfrac{1}{2}\exp(-5 t) + \sin(t)
$$
to repeatedly drive the nominal trajectory toward the safety boundaries. The resulting position trajectories are shown in  Fig.~\ref{fig:states}-Left. The CBF-QP trajectory violates the safety constraint, whereas the trajectory produced by the mirrored Sontag formula~\eqref{eq:sontag_mirror} respect all safety constraints. Fig.~\ref{fig:states}-Right shows the corresponding filtered velocity command $v_{\rm cmd}$. As expected, the proposed smooth formula produces a smooth command, while the CBF-QP command exhibits visible kinks.

To quantify velocity tracking errors, we define $e\triangleq x_2-v_{\rm cmd}$.
The tracking errors and their first two derivatives are shown in~Fig.~\ref{fig:errors}. For the smooth controller, the tracking errors decay to zero after the initial transient. In contrast, for the CBF-QP controller, the nonsmoothness in $v_{\rm cmd}$ repeatedly excites the higher-order tracking dynamics, generating new transients before the previous tracking errors have fully settled. These persistent tracking errors are what ultimately leads to the safety violation observed in Fig.~\ref{fig:states}-Left.

\subsection{Hardware Demonstrations}

We next validate the proposed controller on a Crazyflie~2.1$+$ small quadrotor. In contrast to the deliberately slow tracking dynamics considered in simulation, the lower-level controller of the quadrotor tracks $v_{\rm cmd}$ substantially more closely. The hardware experiments therefore investigate whether the benefits of the proposed smooth selection persist when tracking transients are rapidly attenuated.

For the experiments, the reference $r_\des$ is given as follows. It first smoothly ramps from $0$ to $0.3\, \mathrm{m}$ over $t\in[0,10]$. Then for $t\in[10,35]$, the reference follows the following dynamics:
$$
    \dot r_\des
    =3(1-r_\des)+\sin\!\big(0.3\pi(t-10)\big),
    \qquad r_\des(10)=0.3.
$$
After that, the reference transitions smoothly via a cubic Hermite Interpolation
to $-0.3\,\mathrm{m}$ over $t\in[35,55]$, and for $t\in[55,80]$, it follows the dynamics:
$$
    \dot r_\des
    =3(-1-r_\des)+\sin\!\big(0.3\pi(t-55)\big),
    \qquad r_\des(55)=-0.3.
$$
Finally  $r_\des$ is smoothly returned to zero over $t\in[80,90]$. The two intermediate phases intentionally drive the nominal reference toward and beyond the upper and lower safety boundaries, respectively. This allows us to evaluate the proposed formula~\eqref{eq:sontag_mirror} under each of the two safety constraints.

\begin{figure}
    \centering
    \includegraphics[width=0.9\linewidth]{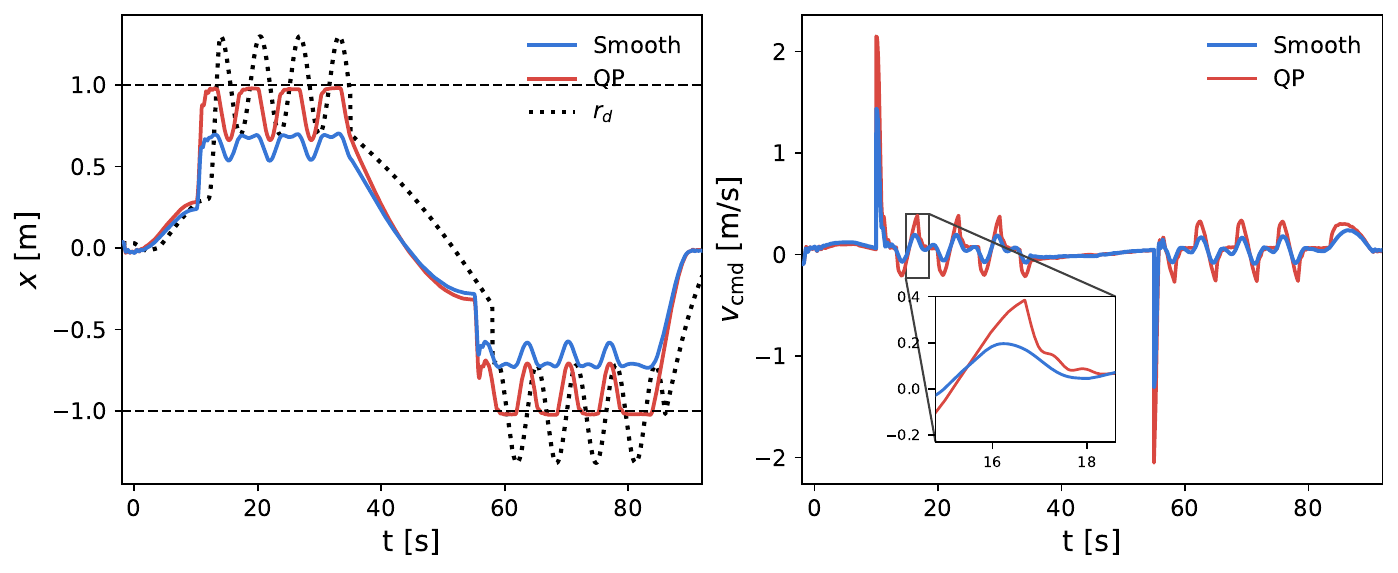}
    \caption{Quadrotor experiment comparing the mirrored Sontag formula~\eqref{eq:sontag_mirror} and the CBF-QP. \textbf{Left:} quadrotor position $x_1$ and reference trajectory $r_d$, with safety boundaries at $x_1=\pm1$.
    \textbf{Right:} corresponding filtered velocity commands $v_{\rm cmd}$. }
    \label{fig:exp_states}
    \vspace{-10pt}
\end{figure}

\begin{figure}
    \centering
    \includegraphics[width=0.9\linewidth]{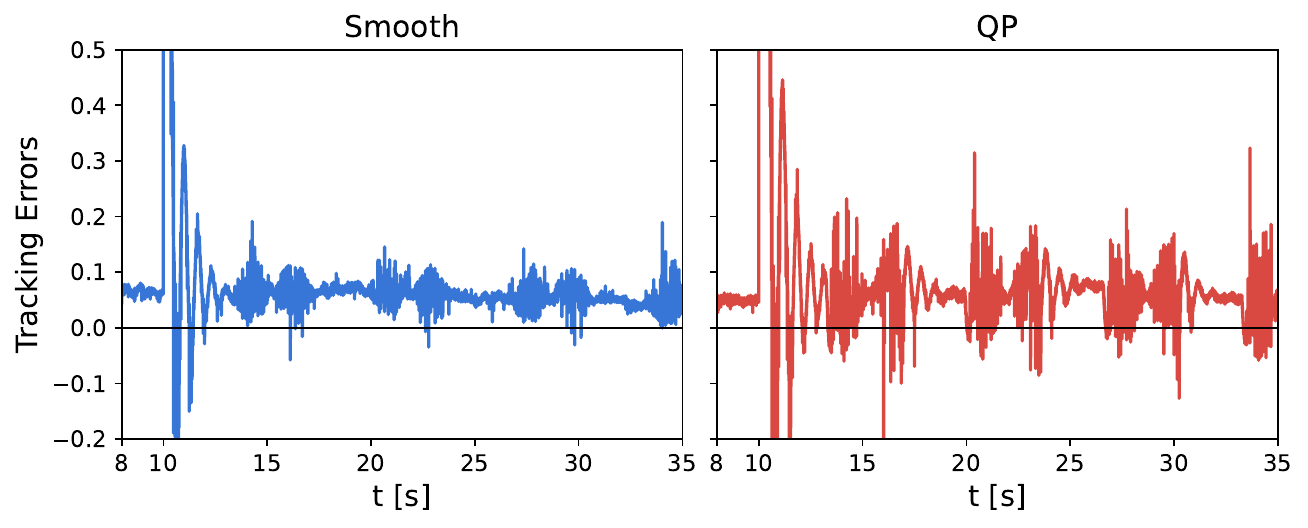}
    \caption{Velocity tracking error $e=x_2-v_{\rm cmd}$ in the hardware experiments under the mirrored Sontag formula~\eqref{eq:sontag_mirror} (\textbf{left}) and the CBF-QP (\textbf{right}) from $t= 8\,\mathrm{s}$ to $t= 35\,\mathrm{s}$.}
    \label{fig:exp_track}
    \vspace{-15pt}
\end{figure}

As shown in Fig.~\ref{fig:exp_states}-Left, the CBF-QP performs substantially better than in simulation and tracks the nominal reference more closely. Nevertheless, its trajectory slightly violates safety constraints, with $x_1\in[-1.025,\,0.983]\,\mathrm{m}$. In contrast, the trajectory generated by the mirrored Sontag formula~\eqref{eq:sontag_mirror} remains within the safe set, with $x_1\in[-0.736,\,0.703] \,\mathrm{m}$.
The small safety violation observed for the CBF-QP is attributed to the mismatch between the single-integrator model used by the safety filter and the actual velocity realized by the quadrotor. Fig.~\ref{fig:exp_track} shows the lower-level controller tracks $v_{\rm cmd}$ closely, but a nonzero tracking error remains due to effects such as measurement noise, localization bias, and other unmodeled dynamics. When the CBF-QP activates the CBF constraints, the QP uses a control on the boundary of the admissible control set, satisfying them without any margin, making it more prone to error.
In contrast, the weighted-centroid construction selects controls from the interior of the admissible set and therefore provides an inherent robustness against such implementation errors, albeit at the cost of more conservative tracking of the nominal reference, as is also evident in Fig.~\ref{fig:exp_states}-Left.

The simulation and hardware results illustrate two benefits of the proposed weighted-centroid construction: smoothness mitigates tracking transients induced by higher-order dynamics, while interior control selection provides additional robustness margin against residual implementation errors.

\vspace*{-1ex}
\section{Conclusion}
We have introduced a weighted-centroid methodology for constructing explicit smooth feedback formulas from state-dependent affine control constraints. The framework recovers Sontag's universal formula as a special case, extends to several multiple-constraint settings, and yields smooth controllers for safety filtering, stabilization, and safe stabilization. Simulation and hardware experiments on a quadrotor illustrate the practical benefits of smooth filtered commands and the inherent robustness margin provided by interior control selection. Future work will investigate broader classes of dependent constraints and more efficient weighted-centroid constructions, and exploit the proposed design in hierarchical control architectures.

{
\bibliographystyle{ieeetr}
\bibliography{
alias,
PO,
main-Pio,
MC}
}

\begin{IEEEbiography}[{\includegraphics[width=1in,height=1.25in,clip,keepaspectratio]{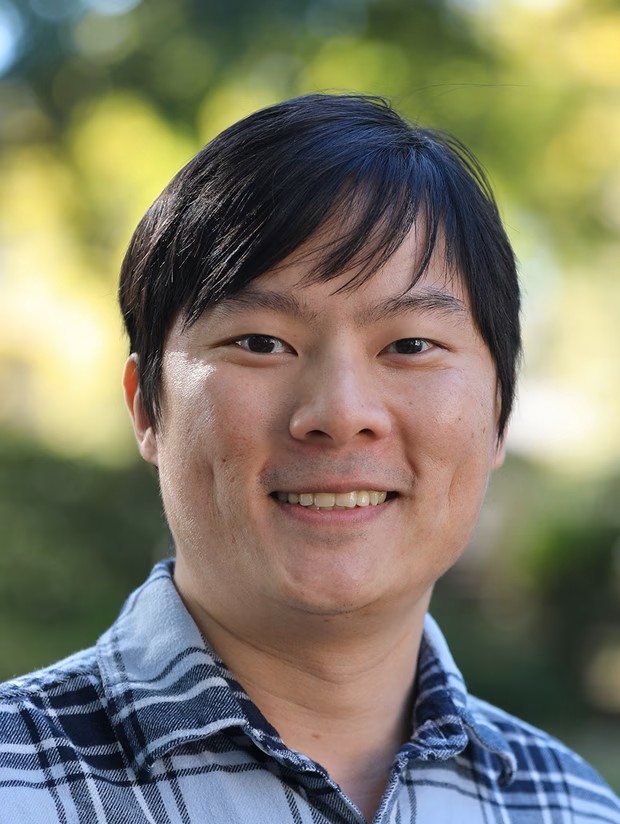}}]{Pio
    Ong} (Member, IEEE)
received the B.S. degree in Aerospace Engineering from University of California, San Diego (UCSD) in 2012, and the M.S. degree in Astronautical Engineering from University of Southern California (USC) in 2013.
He received the Ph.D. degree in Engineering Sciences (Aerospace Engineering) from UCSD in 2022 with Professor Jorge Cort\'{e}s as his advisor. From 2022 to 2026, he was a Postdoctoral Scholar at the California Institute of Technology (Caltech), under the supervision of Professor Aaron D. Ames. He is currently a Postdoctoral Scientist in the Safe Autonomous Foundations (SAF) Lab at Amazon. His current research interests include control barrier functions, safety-critical control, humanoid robot safety, and resilient autonomy.
\end{IEEEbiography}
\vspace*{-7ex}

\begin{IEEEbiography}[{\includegraphics[width=1in,height=1.25in,clip,keepaspectratio]{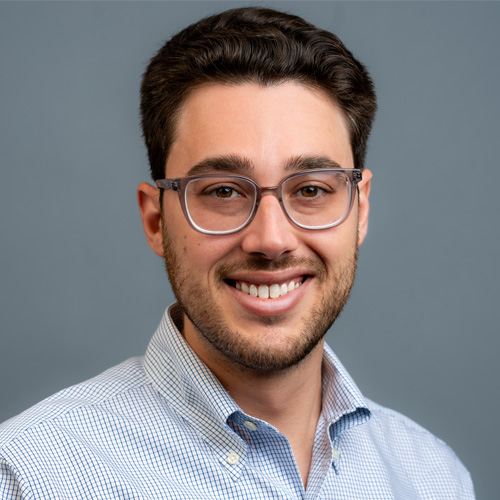}}]{Max H. Cohen} (Member, IEEE)
received the B.S. degree in Mechanical Engineering from the University of Florida in 2018, the M.S. degree in Mechanical Engineering from Boston University in 2022, and the Ph.D. degree in Mechanical Engineering from Boston University in 2023. He was a Postdoctoral Scholar at Caltech from 2023--2025 and is currently an Assistant Professor in the Department of Electrical and Computer Engineering at North Carolina State University, Raleigh, NC.
\end{IEEEbiography}
\vspace*{-7ex}

\begin{IEEEbiography}[{\includegraphics[width=1in,height=1.25in,clip,keepaspectratio]{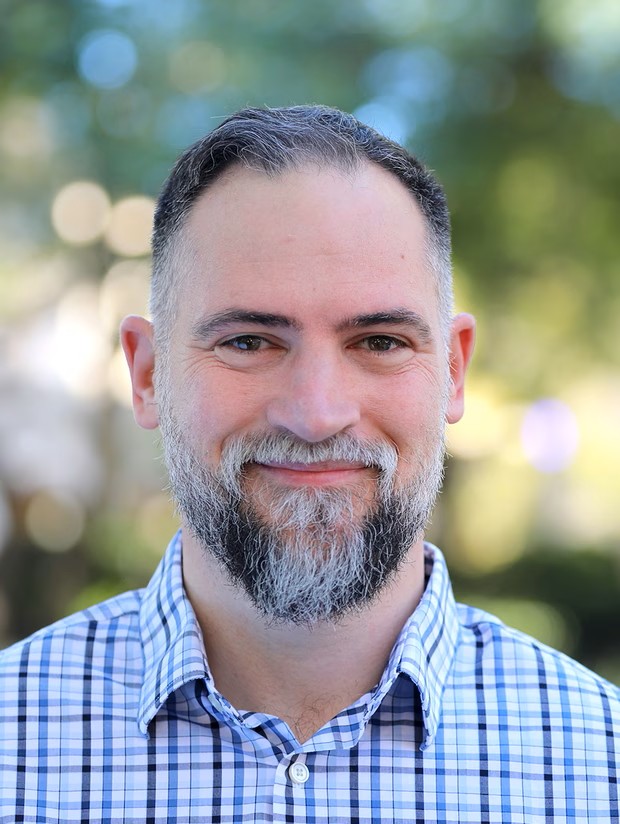}}]{Aaron D. Ames}
(Fellow, IEEE) received the joint B.S. degree in mechanical engineering and B.A. degree in mathematics from the University of St. Thomas in 2001, and the joint M.A. degree in mathematics and the Ph.D. degree in electrical engineering and computer sciences from the University of California, Berkeley, in 2006. He was an Associate Professor with Georgia Tech, Woodruff School of Mechanical Engineering, and the School of Electrical \& Computer Engineering. He is the Bren Professor of mechanical and civil engineering and control and dynamical systems with Caltech. He was a Postdoctoral Scholar in control and dynamical systems with Caltech from 2006 to 2008, and began his faculty career with Texas A\&M University in 2008. His research interests include the areas of robotics, nonlinear, safety-critical control, and hybrid systems, with a special focus on applications to dynamic robots, both formally and through experimental validation.
At UC Berkeley, he was the recipient of the 2005 Leon O. Chua Award for achievement in nonlinear science and the 2006 Bernard Friedman Memorial Prize in Applied Mathematics, NSF CAREER Award in 2010, 2015 Donald P. Eckman Award, and 2019 IEEE CSS Antonio Ruberti Young Researcher Prize.
\end{IEEEbiography}

\vspace*{-7ex}

\begin{IEEEbiography}[{\includegraphics[height=1.1in,clip,keepaspectratio]{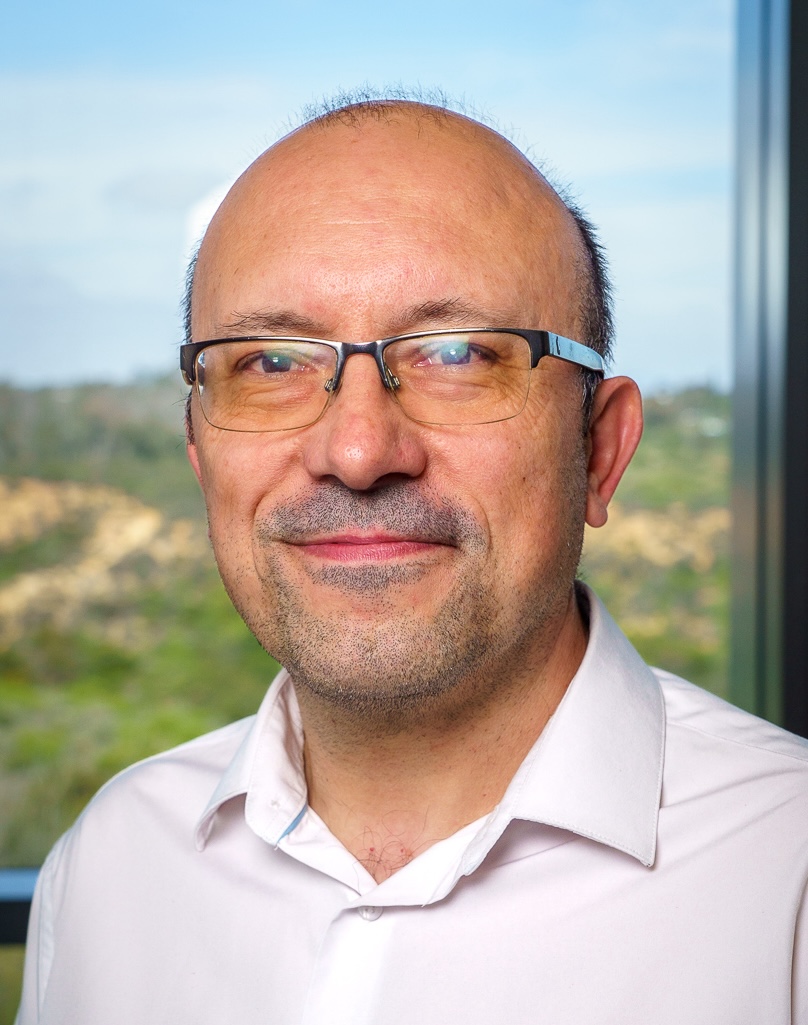}}]{Jorge
    Cort\'{e}s}(M'02, SM'06, F'14) received the Licenciatura degree in
  Mathematics from Universidad de Zaragoza, Spain, in 1997, and the
  Ph.D. degree in Engineering Mathematics from Universidad Carlos III
  de Madrid, Spain, in 2001. He held postdoctoral positions with the
  University of Twente, The Netherlands, and the University of
  Illinois at Urbana-Champaign, USA.
  He is a Professor and Cymer Corporation Endowed Chair in High
  Performance Dynamic Systems Modeling and Control at the Department
  of Mechanical and Aerospace Engineering, UC San Diego, California,
  USA.  He is a Fellow of IEEE, SIAM, and IFAC.  His research
  interests include distributed control and optimization, network
  science, reasoning and decision making under
  uncertainty, and multi-agent coordination in
  robotic, power, and brain networks.
\end{IEEEbiography}

\end{document}